\documentclass[11pt,twoside]{article}
\usepackage[margin=1.3in]{geometry}
\usepackage[T1]{fontenc}
\usepackage{lmodern}
\usepackage{amsmath,amssymb,amsthm,mathtools}
\usepackage{esint}
\usepackage{microtype}
\usepackage{enumitem}
\usepackage{xcolor}
\definecolor{linkblue}{RGB}{25,64,100}
\usepackage[colorlinks=true,linkcolor=linkblue,citecolor=linkblue,urlcolor=linkblue]{hyperref}
\usepackage{bookmark}
\allowdisplaybreaks[2]
\numberwithin{equation}{section}
\newtheorem{theorem}{Theorem}[section]
\newtheorem{lemma}[theorem]{Lemma}
\newtheorem{proposition}[theorem]{Proposition}

\theoremstyle{definition}

\theoremstyle{remark}

\newcommand{\R}{\mathbb R}
\newcommand{\Sn}{\mathbb S^{n-1}}
\newcommand{\eps}{\varepsilon}
\newcommand{\pa}{p^*}
\newcommand{\A}{\mathcal A}
\DeclareMathOperator{\diver}{div}
\DeclareMathOperator{\supp}{supp}

\begin{document} 
\title{ {\bf Asymptotic behavior of solutions to the critical $p$-Laplace equation near isolated singularities} }

\author{Hui Yang}
\date{\today}
\maketitle

\begin{abstract}
In this paper, we study the asymptotic behavior of weak solutions to the following critical $p$-Laplace equation 
\[ 
-\Delta_p u=u^{p^*-1}, \quad u\ge0  \quad\hbox{in }B_1\setminus\{0\}, 
\] 
where $\Delta_pu=\diver(|\nabla u|^{p-2}\nabla u)$ is the classical $p$-Laplacian with $1<p<n$, $p^*=\frac{np}{n-p}$ is the Sobolev critical exponent and $B_1\setminus\{0\}\subset \mathbb{R}^n$ is the punctured unit ball. We develop a nonlinear potential-theoretic method to establish sharp upper bounds for solutions, without using the method of moving planes or moving spheres. Furthermore, we use Steiner symmetrization to study the radial symmetry of global singular solutions, and thereby classify all Fowler singular solutions.  
Finally, we prove that every local singular solution is asymptotic to a global Fowler singular solution. 
\medskip


\end{abstract}

\section{Introduction}
In this paper, we study the asymptotic behavior of weak solutions to the following critical $p$-Laplace equation 
\begin{equation}\label{equ:main26}
-\Delta_p u=u^{p^*-1}, \quad u\ge0  \quad\hbox{in }B_1\setminus\{0\}, 
\end{equation}
where $\Delta_pu=\diver(|\nabla u|^{p-2}\nabla u)$ is the classical $p$-Laplacian with $1<p<n$, $p^*=\frac{np}{n-p}$ is the Sobolev critical exponent and $B_1\setminus\{0\}\subset \mathbb{R}^n$ is the punctured unit ball. The $p$-Laplacian is associated with the $p$-Dirichlet energy $\frac1p\int|\nabla u|^p$, whose first-variation operator is $-\Delta_p$. Its gradient-dependent flux is a basic model of nonlinear diffusion. Related power-growth constitutive laws occur in generalized
viscous-fluid systems, including those introduced in \cite{Ladyzhenskaya1967}. The same energy underlies the theory of minimizing $p$-harmonic maps \cite{HardtLin1987}. The particular source term considered here has a variational origin: after normalization, the critical equation is the Euler--Lagrange equation of the Sobolev quotient
\[
 \frac{\displaystyle\int_{\R^n}|\nabla u|^p\,dx}
 {\left(\displaystyle\int_{\R^n}|u|^{np/(n-p)}\,dx\right)^{(n-p)/n}},
 \qquad 1<p<n.
\]
Its sharp constant and extremals belong to the classical work of Aubin and Talenti \cite{Aubin1976,Talenti1976}. 

For $p=2$, this equation is the Yamabe equation with a flat background metric. More generally,
if $n\ge3$ and $g$ is a Riemannian metric, prescribing positive
constant scalar curvature for $u^{4/(n-2)}g$ leads to
\[
 -\Delta_g u+\frac{n-2}{4(n-1)}R_g u
       =c\,u^{(n+2)/(n-2)},
\]
where $R_g$ is the scalar curvature and $c>0$ is constant.
Caffarelli, Gidas, and Spruck \cite{CGS1989} established the asymptotic behavior of solutions near the singularity in the flat case. Their argument combines measure-theoretic Alexandrov reflection with the radial equation and its asymptotic energy identity. Korevaar, Mazzeo, Pacard, and Schoen \cite{KMPS1999} gave a different
proof and refined the asymptotics using conformally deformed Fowler
profiles, cylinder Jacobi fields, and blow-up analysis. Han, Li, and Li \cite{HanLiLi2021} subsequently established
asymptotic expansions to arbitrary order for positive solutions of the Yamabe equation near isolated singular points. Subsequent to \cite{CGS1989}, other second-order Yamabe type equations with isolated singularities have also been studied; see, e.g., \cite{HLT,LC96,Li06} and the references therein. For smooth background metrics that need not be locally conformally flat, Marques \cite{Marques2008} established Fowler asymptotics in
dimensions $3\le n\le5$, and Xiong and Zhang \cite{XiongZhang2022} treated dimension six. In all of the aforementioned works, the method of moving planes or moving spheres serves as a fundamental tool. The recent papers of Han, Jiang, Wang, Xiong, and Zhang \cite{HJWXZ2026I,HJWXZ2026II} complete this
classification through dimension $29$ and construct counterexamples
for every $n\ge30$. Their Part~I obtains the upper bound from corrected
bubble expansions and an adjoint integral identity, without moving
planes or moving spheres.

For $p\neq2$, the lack of conformal covariance, together with the
singular or degenerate nature of the $p$-Laplacian, makes the classical
moving-plane and moving-sphere methods difficult to adapt to the study
of isolated singularities. Inspired by the approach in Han-Jiang-Wang-Xiong-Zhang \cite{HJWXZ2026I}, we develop a nonlinear potential-theoretic method that does not rely on either technique and establish the sharp upper bound near the singularity, thereby classify the isolated singularities. Throughout the paper, assume that $1<p<n$. 
A nonnegative function $u\in W^{1,p}_{\mathrm{loc}}(B_1\setminus\{0\})$ is called a weak solution of \eqref{equ:main26} if
\[
 \int_{B_1\setminus\{0\}} |\nabla u|^{p-2}\nabla u\cdot\nabla\varphi
 =\int_{B_1\setminus\{0\}} u^{p^*-1} \varphi  \qquad \forall~ \varphi\in C_c^\infty(B_1\setminus\{0\}).
\]

\begin{theorem}[Classification of an isolated singularity]\label{thm:main}
Let $u\in W^{1,p}_{\mathrm{loc}}(B_1\setminus\{0\})$ be a nonnegative weak solution of \eqref{equ:main26}. 
Then there exists a constant $C>0$ such that
\[
u(x)\le C|x|^{-\frac{n-p}{p}} \qquad \text{for all } 0<|x|<1/2. 
\]
Moreover, exactly one of the following alternatives holds.
\begin{enumerate}[label=\textup{(\roman*)},leftmargin=*]
\item The singularity is removable: $u$ has an extension in $W^{1,p}_{\mathrm{loc}}(B_1)\cap C^{1,\alpha}_{\mathrm{loc}}(B_1)$
for some $\alpha\in(0,1)$. 
\item The singularity is nonremovable, and there exist constants $c,C>0$ such that 
\begin{equation}\label{eq:main-two-sided}
 c|x|^{-\frac{n-p}{p}}\le u(x)\le C|x|^{-\frac{n-p}{p}} 
\qquad \text{for all } 0<|x|<1/2. 
\end{equation}
\end{enumerate}
\end{theorem} 

When specialized to the Euclidean metric, the upper-bound argument of Han-Jiang-Wang-Xiong-Zhang \cite{HJWXZ2026I} applies in every dimension $n\ge3$, since the geometric error terms responsible for the restriction 
$n\le29$ vanish. Their Euclidean argument combines an adjoint integral
identity and Green representation with the nondegeneracy of the standard
bubble to obtain an $L^1$ estimate for the bubble error. Our proof of Theorem \ref{thm:main} instead
combines the classification of entire solutions with nonlinear source-mass
estimates, intermediate-scale bubble-tail control, and zero-Pohozaev
rigidity of normalized $p$-harmonic limits, closing the argument through
a second normalization at spherical minima. This nonlinear
potential-theoretic approach establishes the upper bound throughout
$1<p<n$, including both the singular and the degenerate regimes,
without requiring the nondegeneracy of the critical bubble.

To further characterize the asymptotic behavior of singular solutions, we classify the Fowler singular solutions of the critical equation \eqref{equ:main26} on $\R^n\setminus\{0\}$. The key is to establish the radial symmetry of the solutions. Since the classical methods of moving planes \cite{CGS1989} or moving spheres \cite{LiZhang2003,LiZhu1995} are no longer applicable, we use a Steiner symmetrization method to establish the following radial symmetry. 

\begin{theorem}[Radial symmetry of global singular solutions]
\label{thm:global-radial}
Let $U\in W^{1,p}_{\mathrm{loc}}(\R^n\setminus\{0\})$ be a nonnegative weak solution of
\[
 -\Delta_pU=U^{p^*-1} \quad\text{in }\R^n\setminus\{0\}.
\]
If the origin is a nonremovable singularity, then $U(x)=\Phi(|x|)$ with $\Phi'(r)<0$ for every $r>0$, and
\begin{equation}\label{eq:global-automatic-bounds}
 c|x|^{-\frac{n-p}{p}}\le U(x)\le C|x|^{-\frac{n-p}{p}}, \qquad x\ne0 
\end{equation}
for some $c,C>0$.  
\end{theorem} 

For Theorem~\ref{thm:global-radial}, the first step is a rigidity
result under two-sided critical bounds. We localize the continuous
Steiner symmetrization of Brock \cite{Brock2000}, using logarithmic
cutoffs to control the infinite-energy ends. Applying this result
to tangent limits gives a negative Pohozaev constant, which then
forces the global bounds for the original solution. Consequently,
no decay or finite-energy assumption at infinity is needed. Once radial symmetry is established, we can apply the classification of radial solutions due to Franca \cite{Franca2004Classification}. More precisely, put
\[
 a=\frac{n-p}{p}, \qquad v_*=a^a,\qquad E_*=-\frac{a^n}{n},\qquad
 G(v)=-\frac{a^p}{p}v^p+\frac{v^{p^*}}{p^*}.
\]
For each $E\in(E_*,0)$ there is a positive nonconstant periodic
solution $(V_E,Z_E)$ of
\begin{equation}\label{eq:main-fowler-system}
 V_E'=aV_E-Z_E^{1/(p-1)},\qquad Z_E'=V_E^{p^*-1} - aZ_E,
\end{equation}
uniquely normalized by
\[
 V_E(0)=\alpha_E,\qquad Z_E(0)=(a\alpha_E)^{p-1},
 \qquad G(\alpha_E)=E,\quad 0<\alpha_E<v_*.
\]
Write $T(E)>0$ for its minimal period. Then exactly one of the following
representations holds:
\begin{enumerate}[label=\textup{(\roman*)},leftmargin=*]
\item $U(x)=a^a|x|^{-a}$;
\item for a unique $E\in(E_*,0)$ and a phase $\tau$ unique modulo $T(E)$,
\begin{equation}\label{eq:global-fowler-classification}
 U(x)=|x|^{-a}V_E(\log|x|+\tau).
\end{equation}
\end{enumerate}
Conversely, all these functions are positive global solutions with a
nonremovable singularity at the origin. Their Pohozaev constants are
$\omega_{n-1}E_*$ and $\omega_{n-1}E$, respectively.

Building on Theorems \ref{thm:main} and \ref{thm:global-radial}, we can further describe the behavior of solutions near the singularity, which is asymptotically close to a global Fowler singular solution.

\begin{theorem}[Asymptotics to a global singular solution]
\label{thm:asymptotic-radial-symmetry}
Let $u\in W^{1,p}_{\mathrm{loc}}(B_1\setminus\{0\})$ be a nonnegative weak solution of \eqref{equ:main26} with a nonremovable singularity at the origin. 
Then there exists a global singular solution $U_0$ from Theorem~\ref{thm:global-radial} such that
\begin{equation}\label{eq:asymptotic-fixed-profile}
 \lim_{r \to 0}\sup_{\theta\in\Sn}
 \left|\frac{u(r\theta)}{U_0(r)}-1\right|=0.
\end{equation}
\end{theorem}

For Theorem~\ref{thm:asymptotic-radial-symmetry}, subsequential Fowler
limits do not by themselves fix the phase. Our cylinder argument
shares with \cite{KMPS1999} the use of angular Jacobi fields and
Pohozaev invariants, but must control the nonlinear flux explicitly.
We first prove that the rescaled gradients stay away from zero.
Only then do we linearize along secants, obtain angular coercivity,
and derive a quadratic, integrable error in the exact averaged-flux
system. The conserved Pohozaev quantity controls the energy error
at the same order, yielding an integrable phase error and a single
limiting phase. Spherical averaging is therefore not treated as
commuting with $\Delta_p$, and no uniformly elliptic estimate is
used at a vanishing gradient.

We will use standard local $C^{1,\alpha}$ regularity, weak Harnack inequalities for nonnegative $p$-superharmonic functions, and Harnack inequalities for positive solutions of
$-\Delta_pv=V(x)v^{p-1}$ with bounded rescaled potential; see
\cite{DiBenedetto1983,HKM2006,Tolksdorf1984}. The classification of nonnegative entire solutions recently established by Zhang \cite[Theorem~1.1]{Zhang2026} plays an important role in this paper. The classification of entire solutions has developed from the natural
energy space to a fully local weak-solution setting. In the class
$\mathcal D^{1,p}(\R^n)=\{v\in L^{p^*}(\R^n):\nabla v\in L^p(\R^n)\}$,
Damascelli, Merch\'an, Montoro, and Sciunzi \cite{DMMS2014}
established the classification for $1<p<2$ and $p^*\ge2$.
V\'etois \cite{Vetois2016}, using sharp decay estimates and the
moving-plane theorem of Damascelli and Ramaswamy
\cite{DamascelliRamaswamy2001}, removed the restriction $p^*\ge2$;
Sciunzi \cite{Sciunzi2016} treated $2<p<n$, completing the
finite-energy classification. 
Without assuming finite energy, Catino, Monticelli, and Roncoroni
\cite{CatinoMonticelliRoncoroni2023} introduced integral-identity
arguments yielding unconditional classification when $n=2$, $1<p<2$,
and when $n=3$, $3/2<p<2$, as well as results under weaker growth
conditions in other ranges. Ou \cite{Ou2025} obtained classification
without additional global assumptions for $(n+1)/3<p<n$;
V\'etois \cite{Vetois2024} and Sun and Wang \cite{SunWang2025}
further enlarged the admissible exponent ranges. Ciraolo and Gatti \cite{CiraoloGatti2026} classified bounded positive solutions for every $1<p<n$ under the
additional sharp infimum-decay condition
$\inf_{B_R}v\le C R^{-(n-p)/(p-1)}$ for $R\ge1$. Zhang \cite{Zhang2026} proves a complete classification of bounded entire solutions for all $1<p<n$. Together with a 
Schoen-type estimate and the Ciraolo--Gatti classification,  Zhang \cite{Zhang2026} further obtains a fully unrestricted classification theorem.

The paper is organized as follows. In Section \ref{sec:upper}, we establish upper bounds for solutions near the singularity.  
In Section \ref{sec:classification-proof}, we classify the isolated singularities and thereby prove Theorem \ref{thm:main}. In Section \ref{sec:global-symmetry}, we show the radial symmetry of global solutions under two-sided critical bounds. In Section \ref{sec:global-without-bounds}, we establish the radial symmetry and the sharp estimates in Theorem \ref{thm:global-radial}. In Section \ref{sec:asymptotic-symmetry}, we prove Theorem \ref{thm:asymptotic-radial-symmetry}.

\section{The sharp upper bound from the entire classification}
\label{sec:upper}

In this section the only classification result for the critical nonlinear
equation that we use is Zhang's Theorem~1.1 in \cite{Zhang2026}. Ordinary Harnack inequalities are used only after the rescaled potential has been bounded.
Set 
\[
 a=\frac{n-p}{p},\qquad b=\frac{n-p}{p-1},\qquad
 q=\pa-1,\qquad p'=\frac p{p-1},\qquad
 \sigma=q-p+1=\frac pa.
\]
We shall repeatedly use
\begin{equation}\label{eq:ub-exponents}
 b=ap',\qquad b(p-1)=n-p,\qquad bq=n+p',\qquad
 aq=n-a.
\end{equation}
Write $\omega_{n-1}=|\Sn|$. We denote the normalized bubble by
\begin{equation}\label{eq:entire-bubble}
 \mathcal U(x)=(1+\gamma_{n,p}|x|^{p'})^{-a},
 \qquad \gamma_{n,p}=\frac{1}{b\,n^{1/(p-1)}}.
\end{equation}
Put
\[
 \kappa=\lim_{r\to\infty}r^b\mathcal U(r)>0.
\]
Its explicit formula and its radial equation give
\begin{equation}\label{eq:ub-bubble-mass}
 \mathcal U(r)\le\kappa r^{-b},\qquad
 \int_{\R^n}\mathcal U^q\,dx
   =\omega_{n-1}(b\kappa)^{p-1}.
\end{equation}
Notice that the integrability in this formula concerns $\mathcal U^q$;
it follows from $bq=n+p'>n$.

\begin{theorem}\label{thm:upper}
Let $1<p<n$ and let $u\in W^{1,p}_{\mathrm{loc}}(B_1\setminus\{0\})$ be a nonnegative weak solution of \eqref{equ:main26}. Then there exists $C=C(u,n,p)>0$ such that
\[
 u(x)\le C|x|^{-\frac{n-p}{p}},\qquad 0<|x|<\tfrac12.
\]
\end{theorem}

We first establish the estimates and compactness statements needed in
the proof. We use the usual lower semicontinuous representatives of
$p$-superharmonic functions. Their generalized gradients and Riesz
measures are understood in the standard truncation sense; see
\cite{HKM2006,KKT2011}. A nonnegative $p$-superharmonic function on the
complement of a set of zero $p$-capacity extends across that set as a
$p$-superharmonic function, provided it is not identically infinite.
We use this extension theorem in the form recorded in
\cite[Theorem~6.3(d)]{Bjorn2006}. 

\subsection{A source-mass estimate} 
\begin{lemma}\label{lem:ub-mass}
Suppose that $v\ge0$ is $p$-superharmonic in $B_{8r}(z)$, with Riesz
measure $\mu=-\Delta_pv$. Then
\begin{equation}\label{eq:ub-mass}
 \mu(B_r(z))\le C r^{n-p}
       \left(\inf_{B_{2r}(z)}v\right)^{p-1},
\end{equation}
where $C$ depends only on $n$ and $p$. In particular, if $v(z)\le M$,
then $\mu(B_r(z))\le C M^{p-1}r^{n-p}$.
\end{lemma}

\begin{proof}
First suppose that $v$ is a bounded weak supersolution. Choose
\[
 0<\tau<\min\left\{p-1,\frac p{n-p}\right\}
\]
and a cutoff supported in $B_{4r}(z)$, equal to one on $B_{2r}(z)$,
with gradient bounded by $C/r$. Testing with
$\eta^p(v+\delta)^{-\tau}$ and using Young's inequality gives
\[
 \int_{B_{2r}(z)}(v+\delta)^{-1-\tau}|\nabla v|^p
 \le C r^{-p}\int_{B_{4r}(z)}(v+\delta)^{p-1-\tau}.
\]
The weak Harnack inequality for nonnegative supersolutions is valid for
every exponent strictly below $n(p-1)/(n-p)$. Applying it on these
nested balls and then letting $\delta\downarrow0$, we obtain, with
$m=\inf_{B_{2r}(z)}v$,
\[
 \int_{B_{2r}(z)}v^{-1-\tau}|\nabla v|^p
       \le C r^{n-p}m^{p-1-\tau}.
\]
The same calculation can throughout be kept in its $v+\delta$ form
if $m=0$. Since
$(1+\tau)(p-1)<n(p-1)/(n-p)$, H\"older's inequality and a second
application of weak Harnack yield
\begin{align}
 \int_{B_{2r}(z)}|\nabla v|^{p-1}
 &\le
 \left(\int_{B_{2r}(z)}v^{-1-\tau}|\nabla v|^p\right)^{(p-1)/p}
 \left(\int_{B_{2r}(z)}v^{(1+\tau)(p-1)}\right)^{1/p}
 \notag\\
 &\le C r^{n-p+1}m^{p-1}.
 \label{eq:ub-flux-L1}
\end{align}
Testing the Riesz measure with a nonnegative cutoff supported in
$B_{2r}(z)$ and equal to one on $B_r(z)$ proves
\eqref{eq:ub-mass} in this case.

For a general $p$-superharmonic function, let $v_k=\min\{v,k\}$.
These are locally Sobolev weak supersolutions. If $Dv$ denotes the
generalized gradient, then
\[
 |\nabla v_k|^{p-1}=|Dv|^{p-1}\mathbf1_{\{v<k\}}
 \quad\text{a.e.}
\]
The right side of \eqref{eq:ub-flux-L1}, applied to $v_k$, is bounded
by the expression with $m=\inf_{B_{2r}(z)}v$. Monotone convergence
therefore gives $|Dv|^{p-1}\in L^1$ on smaller balls. Dominated
convergence then gives
\[
 |\nabla v_k|^{p-2}\nabla v_k\longrightarrow
 |Dv|^{p-2}Dv\quad\text{in }L^1_{\mathrm{loc}}.
\]
Consequently the Riesz measures of $v_k$ converge weakly to that of
$v$. The lower semicontinuity of the mass of an open ball proves
\eqref{eq:ub-mass}. The last assertion follows from
$\inf_{B_{2r}(z)}v\le v(z)$.
\end{proof}

\subsection{Pohozaev flux and the constant term at a harmonic pole}

\begin{lemma}\label{lem:pohozaev}
Let $w\ge0$ be a locally $C^{1,\alpha}$ weak solution of
$-\Delta_pw=\eps w^q$, with $\eps\ge0$, in an annulus centered at
$z$. Define
\begin{align}
 \mathcal P_\eps(w,r;z)=\int_{\partial B_r(z)}
 \bigg[&a w|\nabla w|^{p-2}\partial_\nu w\notag\\
 &+r\left\{|\nabla w|^{p-2}(\partial_\nu w)^2
       -\frac{|\nabla w|^p}{p}
       +\frac{\eps w^{\pa}}{\pa}\right\}\bigg]\,dS.
 \label{eq:ub-pohozaev}
\end{align}
Expressions containing a vanishing gradient are interpreted by
continuous extension. This flux is independent of $r$. If $w$ is a
solution in a ball centered at $z$, its flux is zero.
\end{lemma}

\begin{proof}
We recall a weak justification to avoid assuming second derivatives
at critical points of $w$. Freeze $f(x)=\eps w(x)^q$. The weak equation
and convexity imply that $w$ minimizes
\[
 v\longmapsto\int\left(\frac{|\nabla v|^p}{p}-f(x)v\right)
\]
under compactly supported Sobolev perturbations. Its inner-variation
identity is therefore valid. More precisely, for a compactly supported
$C^1$ vector field $X$, use the competitors
$w\circ(\operatorname{Id}+tX)^{-1}$. Changing variables in the gradient
term and using the $C^1$ regularity of $w$ in the frozen source term gives
\[
 \int\left[\left(\frac{|\nabla w|^p}{p}
                  -\frac{\eps w^{\pa}}{\pa}\right)\diver X
 -|\nabla w|^{p-2}\nabla w\cdot DX\nabla w\right]=0.
\]
Here the source contribution was integrated using
$f\nabla w=\nabla(\eps w^{\pa}/\pa)$.
Take $X=(x-z)\eta$ and combine this identity with the weak equation
tested against $aw\eta$, for a scalar compactly supported cutoff
$\eta$. Since $n/\pa=a$, this gives, in distributions,
\[
 \diver\left[
 ((x-z)\cdot\nabla w+aw)|\nabla w|^{p-2}\nabla w
 -(x-z)\left(\frac{|\nabla w|^p}{p}
                    -\frac{\eps w^{\pa}}{\pa}\right)
 \right]=0.
\]
Radial cutoffs give conservation for almost every radius. The local
$C^1$ regularity makes the boundary flux continuous in the radius,
so conservation holds for every admissible radius. In a full ball,
boundedness of $w$ and $\nabla w$ near its center makes the flux tend
to zero as the radius tends to zero.
\end{proof}

\begin{lemma}\label{lem:ub-pole}
Suppose that $W>0$ is $p$-harmonic in $B_R\setminus\{0\}$ and has a
nonremovable singularity at zero. Then there are $k>0$ and $B\in\R$
such that
\begin{equation}\label{eq:ub-pole-expansion}
 W(x)=k|x|^{-b}+B+o(1),\qquad
 |x|\left|\nabla(W-k|x|^{-b})\right|\longrightarrow0.
\end{equation}
Moreover,
\begin{equation}\label{eq:ub-pole-flux}
 \mathcal P_0(W,r;0)=-a\omega_{n-1}(bk)^{p-1}B.
\end{equation}
If $B=0$ and $W\ge k|x|^{-b}$ in the punctured ball, then
$W=k|x|^{-b}$ there.
\end{lemma}

\begin{proof}
The classical isolated-singularity theorem gives
\begin{equation}\label{eq:ub-classical-pole}
 W=k|x|^{-b}+g,\qquad g\in L^\infty(B_R),\qquad
 |x|^{b+1}\nabla g(x)\longrightarrow0;
\end{equation}
see Kichenassamy--V\'eron \cite{KV1986}, as recorded in
\cite[Theorem~2.1]{Veron2002}. The assertion about the constant
term in \eqref{eq:ub-pole-expansion} is stronger than merely bounded
remainder, and we prove it here.

Let $M_g(r)=\max_{\partial B_r}g$ and
$m_g(r)=\min_{\partial B_r}g$. Comparison on annuli with
$(k\pm\delta)|x|^{-b}+M_g(R_1)$ and
$(k\pm\delta)|x|^{-b}+m_g(R_1)$, where $R_1<R$, gives
\[
 M_g(r)\le M_g(R_1),\qquad m_g(r)\ge m_g(R_1),
 \qquad 0<r<R_1.
\]
Indeed, the appropriate $\delta|x|^{-b}$ term dominates the bounded
remainder on a sufficiently small inner sphere; first let that
sphere shrink to zero and then let $\delta\downarrow0$.
Thus $M_g$ and $m_g$ have finite inward limits $B_+\ge B_-$.

Set $\Phi(y)=k|y|^{-b}$,
$U_r(y)=r^bW(ry)$, and $g_r(y)=g(ry)$. On every fixed compact
annulus, $U_r=\Phi+r^bg_r$, and subtraction of the two harmonic
equations gives
\[
 \diver(A_r\nabla g_r)=0,\qquad
 A_r=\int_0^1D\mathcal A\bigl(\nabla\Phi
               +t(\nabla U_r-\nabla\Phi)\bigr)\,dt,
 \qquad \mathcal A(\xi)=|\xi|^{p-2}\xi.
\]
By \eqref{eq:ub-classical-pole}, $\nabla U_r\to\nabla\Phi$ uniformly
on compact annuli, and $\nabla\Phi$ never vanishes there. Local
$p$-harmonic regularity gives uniform $C^{1,\alpha}$ bounds for $U_r$.
Consequently $A_r$ are uniformly elliptic and uniformly H\"older
continuous on each such annulus, also when $p<2$.

The functions $g_r$ are uniformly bounded. Linear interior estimates
give subsequential $C^1$ compactness on compact annuli. Every limit
has maximum $B_+$ and minimum $B_-$ on every centered sphere, and
lies between these same constants everywhere. The strong maximum
principle forces it to be constant. Hence $B_+=B_-=B$ and $g\to B$.
Applying the same linear gradient estimate to $g_r-B$ gives
\[
 \|\nabla g_r\|_{L^\infty(K')}
       \le C\|g_r-B\|_{L^\infty(K)}\longrightarrow0
\]
on nested compact annuli. This proves \eqref{eq:ub-pole-expansion}.

Substitution into \eqref{eq:ub-pohozaev} now proves
\eqref{eq:ub-pole-flux}. More explicitly, the pure fundamental-solution
terms cancel because $(p-1)b/p=a$. The term containing $B$ is
$-a\omega_{n-1}(bk)^{p-1}B$. All terms containing the remainder of
$g$ are $o(1)$: the gradient remainder is $o(r^{-1})$, and
$(b+1)(p-1)=n-1$ makes its products with the leading gradient contribute
only $o(1)$ after integration over the sphere. Conservation then gives
the stated identity at every radius.

Finally suppose that $B=0$ and $W\ge\Phi$. At a contact point,
$\nabla W=\nabla\Phi\ne0$. In a neighborhood of that point the
difference equation can be linearized with uniformly elliptic
coefficients. The strong maximum principle shows that the contact
set is open as well as closed. If there is no contact point, then
$W>\Phi$. On a fixed smaller sphere there is a positive gap
$\delta_0$. Comparison with
$(k-\delta)|x|^{-b}+\delta_0$, first shrinking the inner sphere and
then sending $\delta\downarrow0$, gives $W-\Phi\ge\delta_0$ inside
that sphere. This contradicts $g\to0$.
\end{proof}

\subsection{A bubble-tail estimate at intermediate scales}

\begin{lemma}\label{lem:ub-tail}
Suppose $\Gamma_j\to\infty$ and
\[
 -\Delta_pv_j=v_j^q\quad\text{in }B_{4\Gamma_j},\qquad
 v_j(0)=1,\qquad 0<v_j\le C_0,
 \qquad v_j\to\mathcal U\quad\text{in }C^1_{\mathrm{loc}}(\R^n).
\]
For every fixed $\theta>1$ and every sequence
$S_j\to\infty$ with $S_j/\Gamma_j\to0$, one has
\begin{equation}\label{eq:ub-tail}
 v_j<\theta\mathcal U\quad\text{in }B_{S_j}
\end{equation}
for all sufficiently large $j$.
\end{lemma}

\begin{proof}
Otherwise there are first-contact radii $\rho_j\le S_j$ and
$e_j\in\Sn$ such that
\[
 v_j<\theta\mathcal U\text{ in }B_{\rho_j},\qquad
 v_j(\rho_je_j)=\theta\mathcal U(\rho_j).
\]
Local convergence implies $\rho_j\to\infty$, whereas
$\rho_j/\Gamma_j\to0$. Put
\[
 W_j(x)=\rho_j^bv_j(\rho_jx),\qquad
 \eps_j=\rho_j^{-p'},\qquad \nu_j=\eps_jW_j^q\,dx.
\]
Then
\begin{equation}\label{eq:ub-tail-scaled}
 -\Delta_pW_j=\eps_jW_j^q,\qquad
 \eps_j^aW_j\le C_0,\qquad
 W_j(e_j)\longrightarrow\theta\kappa.
\end{equation}
Inside $B_1\setminus\{0\}$ we have
$W_j\le\theta\kappa|x|^{-b}$. Interior regularity therefore gives,
after passage to a subsequence,
\[
 W_j\to W\quad\text{in }C^1_{\mathrm{loc}}(B_1\setminus\{0\}),
 \qquad -\Delta_pW=0.
\]
For every compactly supported continuous test function in $B_1$,
the change of variables $y=\rho_jx$ gives
\[
 \int\varphi\,d\nu_j
 =\int_{B_{\rho_j}}\varphi(y/\rho_j)v_j(y)^q\,dy
 \longrightarrow
 \varphi(0)\int_{\R^n}\mathcal U^q\,dy.
\]
Here the bound $v_j^q\le\theta^q\mathcal U^q$ on $B_{\rho_j}$ and
\eqref{eq:ub-bubble-mass} justify dominated convergence. The same
argument with the characteristic function of a fixed smaller ball,
together with boundary $C^1$ convergence, gives
\[
 -\int_{\partial B_r}|\nabla W|^{p-2}\partial_\nu W\,dS
       =\omega_{n-1}(b\kappa)^{p-1},\qquad 0<r<1.
\]
Thus $W$ is nonzero, positive, and has a genuine harmonic pole with
coefficient $\kappa$. Since every $W_j$ solves its critical equation
in a full ball, Lemma~\ref{lem:pohozaev} and boundary convergence
give $\mathcal P_0(W,r;0)=0$. Lemma~\ref{lem:ub-pole} shows that its constant
remainder is zero.

We also need a lower comparison. Set $L_j=\Gamma_j/\rho_j\to\infty$.
Given $0<\delta<\kappa$, choose
$r_\delta<\min\{\delta,1/2\}$ so small that
\[
 W\ge(\kappa-\delta/2)r_\delta^{-b}
 \quad\text{on }\partial B_{r_\delta}.
\]
For large $j$, comparison of the supersolution $W_j$ with the
$p$-harmonic function
\[
 (\kappa-\delta)(|x|^{-b}-L_j^{-b})
\]
is valid on $B_{L_j}\setminus\overline B_{r_\delta}$: the inner
boundary inequality follows from convergence, and the outer one
from positivity. Sending $j\to\infty$ and then $\delta\downarrow0$
shows that $W\ge\kappa|x|^{-b}$. The last assertion of
Lemma~\ref{lem:ub-pole} yields
\begin{equation}\label{eq:ub-tail-inner-limit}
 W(x)=\kappa|x|^{-b}\quad(0<|x|<1).
\end{equation}

It remains to justify convergence at the moving contact points;
the interior convergence just established does not suffice for this.
By \eqref{eq:ub-tail-scaled}, $W_j(e_j)\le M$ uniformly. The
mass estimate gives, for every sufficiently small fixed $r>0$,
\begin{equation}\label{eq:ub-contact-mass}
 \nu_j(B_r(e_j))\le C M^{p-1}r^{n-p}.
\end{equation}
Define the intrinsic inverse length
\[
 \Theta_j=\eps_j^{1/p}W_j^{1/a}.
\]
The following elementary core estimate will be used twice. At a point
$z_j$, let $H_j=W_j(z_j)$ and $\ell_j=\Theta_j(z_j)^{-1}$.
If $\Theta_j\le2\Theta_j(z_j)$ in $B_{2\ell_j}(z_j)$, then
\[
 Z_j(y)=H_j^{-1}W_j(z_j+\ell_jy)
\]
solves $-\Delta_pZ_j=Z_j^q$ in $B_2$, has $Z_j(0)=1$, and is
bounded above by $2^a$. Uniform interior regularity makes
$Z_j\ge1/2$ on a fixed ball $B_c$. Consequently
\begin{equation}\label{eq:ub-core-mass}
 \nu_j(B_{c\ell_j}(z_j))
 \ge c_0\eps_jH_j^q\ell_j^n
 =\frac{c_0}{\eps_j^aH_j}
 \ge\frac{c_0}{C_0}=:m_*>0.
\end{equation}

Choose $r$ so small that the right side of
\eqref{eq:ub-contact-mass} is less than $m_*$. If $\Theta_j$ were
unbounded in $B_{r/2}(e_j)$, maximize
$\operatorname{dist}(x,\partial B_r(e_j))\Theta_j(x)$ in
$B_r(e_j)$. At a maximizing point $z_j$, with boundary distance
$d_j$, one has $d_j\Theta_j(z_j)\to\infty$.
The ball $B_{2\ell_j}(z_j)$ is then contained in $B_r(e_j)$ and
satisfies the doubling assumption of \eqref{eq:ub-core-mass}.
This contradicts \eqref{eq:ub-contact-mass}.
Thus $\Theta_j$ is uniformly bounded in a fixed contact neighborhood.
Only now do we apply the ordinary bounded-potential Harnack inequality
to
\[
 -\Delta_pW_j=\Theta_j^pW_j^{p-1}.
\]
Together with $W_j(e_j)\le M$, it gives uniform local bounds and then
$C^1$ compactness there. Taking $e_j\to e\in\Sn$, the local limit
agrees with \eqref{eq:ub-tail-inner-limit} on the inner side of the
unit sphere. Continuity therefore gives $W_j(e_j)\to\kappa$, contrary
to \eqref{eq:ub-tail-scaled}.
\end{proof}

\subsection{A large value forces a small spherical minimum}

Assume throughout the rest of this section that $u\not\equiv0$.
The strong minimum principle and interior regularity give $u>0$ and
$u\in C^{1,\alpha}_{\mathrm{loc}}(B_1\setminus\{0\})$; see
\cite{DiBenedetto1983,Tolksdorf1984,Vazquez1984}. Extend $u$
$p$-superharmonically across zero, and define
\[
 m(r)=\min_{|x|=r}u(x),\qquad h(r)=r^am(r).
\]
Both functions are positive and continuous. The minimum principle
inside $B_r$, followed by comparison with a radial harmonic function
outside $B_r$, shows that for every fixed $R\ge2$
\begin{equation}\label{eq:ub-normalized-lower}
 \frac{u(ry)}{m(r)}\ge c_R>0\qquad (|y|\le R)
\end{equation}
when $r$ is sufficiently small. For clarity, the lower bound on
$|y|\le1$ is one. On $1<|y|<2R$ use
\[
 \frac{|y|^{-b}-(2R)^{-b}}{1-(2R)^{-b}},
\]
whose minimum on $1\le|y|\le R$ is positive.

Suppose that there are $x_j\to0$ such that
\begin{equation}\label{eq:ub-bad-points}
 r_j=|x_j|,\qquad K_j=r_j^au(x_j)\longrightarrow\infty.
\end{equation}
In $D_j=B_{r_j/4}(x_j)$ maximize
$d(x)^{n-p}u(x)$, where $d(x)=\operatorname{dist}(x,\partial D_j)$.
Write $\xi_j$ for a maximizer, $H_j=u(\xi_j)$, and $d_j=d(\xi_j)$.
Then
\begin{equation}\label{eq:ub-selected-product}
 H_j\ge u(x_j),\qquad
 H_jd_j^{n-p}\ge c r_j^{a(p-1)}K_j.
\end{equation}
Set
\[
 \delta_j=H_j^{-1/a},\qquad \Gamma_j=\frac{d_j}{8\delta_j},
 \qquad v_j(y)=\frac{u(\xi_j+\delta_jy)}{H_j}.
\]
The selection gives $v_j\le2^{n-p}$ on $B_{4\Gamma_j}$, and
\begin{equation}\label{eq:ub-selected-radius}
 \Gamma_j^{n-p}=c d_j^{n-p}H_j^p\ge cK_j^p\longrightarrow\infty.
\end{equation}
Moreover, on each fixed compact set,
\[
 v_j(y)\le
 \left(\frac{d_j}{d_j-\delta_j|y|}\right)^{n-p}\longrightarrow1.
\]
Thus every convergent subsequence has a positive entire limit bounded
by one and equal to one at the origin. Zhang \cite[Theorem~1.1]{Zhang2026}
identifies it with $\mathcal U$, so, after extraction,
$v_j\to\mathcal U$ in $C^1_{\mathrm{loc}}$.

Choose
\[
 \lambda_j=K_j^{-1/[2(n-p)]},\qquad S_j=\lambda_j\Gamma_j.
\]
Then $S_j\to\infty$ by \eqref{eq:ub-selected-radius}, and
$S_j=o(\Gamma_j)$. Lemma~\ref{lem:ub-tail} gives, on the 
sphere of radius $\delta_jS_j=d_j\lambda_j/8$ about $\xi_j$,
\[
 u\le CH_jS_j^{-b}
   =C\lambda_j^{-b}(H_jd_j^{n-p})^{-1/(p-1)}.
\]
This sphere is contained in $D_j$. The lower bound
\eqref{eq:ub-normalized-lower}, applied with scale $r_j$, and
\eqref{eq:ub-selected-product} therefore imply
\begin{equation}\label{eq:ub-small-minimum}
 h(r_j)\le CK_j^{-1/[2(p-1)]}\longrightarrow0.
\end{equation}

\subsection{Compactness at a small spherical minimum}

Consider any sequence $t_j\downarrow0$ such that
$\eta_j=h(t_j)\to0$. Define
\begin{equation}\label{eq:ub-min-normalization}
 w_j(y)=\frac{u(t_jy)}{m(t_j)},\qquad
 \eps_j=\eta_j^{p/a}\longrightarrow0.
\end{equation}
Then, away from zero in expanding balls,
\[
 -\Delta_pw_j=\eps_jw_j^q,\qquad \min_{\partial B_1}w_j=1.
\]
Let $\mu_j$ be the Riesz measure of the superharmonic extension of
$w_j$; thus $\mu_j=\eps_jw_j^q\,dy$ off zero. Fix once and for all
$0<\gamma<n(p-1)/(n-p)$. For every fixed $R\ge2$, weak Harnack,
Lemma~\ref{lem:ub-mass}, and \eqref{eq:ub-normalized-lower} give
\begin{equation}\label{eq:ub-basic-compactness}
 w_j\ge c_R>0\text{ on }B_R,\qquad
 \int_{B_R}w_j^\gamma\le C_R,\qquad
 \mu_j(B_R)\le C_R.
\end{equation}
For the last two estimates use $\inf_{B_{2R}}w_j\le1$.
The measure bound includes any atom at zero.

We first claim that
\begin{equation}\label{eq:ub-amplitude-cap}
 \sup_j\sup_K\eps_j^aw_j<\infty, 
 \qquad K\Subset\R^n\setminus\{0\}.
\end{equation}
If this failed, take a fixed ball $D$ compactly contained in
$\R^n\setminus\{0\}$, containing violating points a fixed positive
distance from its boundary. Maximize $d^{n-p}w_j$ in $D$, and denote
the selected height and boundary distance by $H_j,d_j$. Then
\[
 Q_j:=\eps_j^aH_jd_j^{n-p}\longrightarrow\infty.
\]
The intrinsic scale and available radius are
\[
 \delta_j=\eps_j^{-1/p}H_j^{-1/a},\qquad
 \Gamma_j=\frac{d_j}{8\delta_j}.
\]
The selected normalized functions solve the pure critical equation,
are bounded by $2^{n-p}$ on $B_{4\Gamma_j}$, and converge to
$\mathcal U$ by Zhang \cite[Theorem~1.1]{Zhang2026}. Indeed,
\[
 \Gamma_j^{n-p}=cQ_jH_j^{p-1}\longrightarrow\infty,
\]
since \eqref{eq:ub-basic-compactness} supplies a positive lower bound
for $H_j$. Set $S_j=\Gamma_jQ_j^{-1/[2(n-p)]}$. Then
$S_j\to\infty$, $S_j=o(\Gamma_j)$, and the bubble-tail estimate,
on a physical sphere still contained in $D$, gives
\[
 w_j\le CH_jS_j^{-b}=CQ_j^{-1/[2(p-1)]}\longrightarrow0.
\]
This contradicts the positive lower bound in
\eqref{eq:ub-basic-compactness}, and proves
\eqref{eq:ub-amplitude-cap}.

After a diagonal extraction let $\mu_j\rightharpoonup\mu$ locally,
and write $\Theta_j=\eps_j^{1/p}w_j^{1/a}$. Define
\begin{equation}\label{eq:ub-concentration-set}
 \Sigma=\left\{z\ne0:
  \text{for every }0<\rho<|z|/2,\quad
  \limsup_{j\to\infty}\sup_{B_\rho(z)}\Theta_j=\infty
          \right\}.
\end{equation}
This set is relatively closed away from zero. If
$x_j\to z\ne0$ and $\Theta_j(x_j)\to\infty$, perform distance-weighted
selection in $B_{\Theta_j(x_j)^{-1/2}}(x_j)$.
The selected intrinsic core has center tending to $z$, radius tending
to zero, and satisfies the doubling condition used in
\eqref{eq:ub-core-mass}. That calculation and
\eqref{eq:ub-amplitude-cap} give
\[
 \mu_j(B_{c\ell_j}(z_j))
       \ge\frac{c_0}{\eps_j^aH_j}\ge c_K>0
\]
when the centers remain in a fixed compact annulus $K$.
Testing the measures against continuous functions equal to one near
$z$, and then shrinking their supports, proves
$\mu(\{z\})\ge c_K$. A diagonal selection of such $x_j$ applies to
each $z\in\Sigma$. The same constant works locally uniformly in
$z$, so \eqref{eq:ub-basic-compactness} implies that $\Sigma$ is finite
in every compact annulus.

Outside $\Sigma\cup\{0\}$, definition
\eqref{eq:ub-concentration-set} gives a local uniform bound for
$\Theta_j$. The bounded-potential Harnack inequality for
$-\Delta_pw_j=\Theta_j^pw_j^{p-1}$, together with the integral bound
in \eqref{eq:ub-basic-compactness}, now gives local uniform bounds for
$w_j$. Hence, after further extraction,
\begin{equation}\label{eq:ub-harmonic-limit}
 w_j\to W\quad\text{in }
 C^1_{\mathrm{loc}}(\R^n\setminus(\Sigma\cup\{0\})),\qquad
 -\Delta_pW=0,\qquad W>0.
\end{equation}
The measure $\mu$ is supported on $\Sigma\cup\{0\}$, because
$\eps_jw_j^q\to0$ uniformly on every regular compact set.
For $z\in\Sigma$, choose a small ball containing no other exceptional
point. Its boundary has zero $\mu$-mass. Weak convergence of measures
and boundary $C^1$ convergence give
\begin{equation}\label{eq:ub-atom-flux}
 \mu(\{z\})
 =-\int_{\partial B_\rho(z)}
       |\nabla W|^{p-2}\partial_\nu W\,dS>0.
\end{equation}
Consequently each point of $\Sigma$ is a genuine isolated harmonic
pole of $W$.

We will also need the following uniform convergence of spherical
minima:
\begin{equation}\label{eq:ub-minima-convergence}
 \min_{|y|=s}w_j(y)\longrightarrow\min_{|y|=s}W(y)
 \quad\text{uniformly on every compact interval in }(0,\infty).
\end{equation}
At poles the value of $W$ in this formula is understood to be infinity.
To prove it, enlarge the compact annulus under consideration slightly.
It contains only finitely many poles. Choose one radial segment
avoiding all of them; convergence on this segment gives a common
upper bound $M$ for all spherical minima. Around each pole choose
a small ball, avoiding the others and zero, on whose boundary
$W>M+2$. Boundary convergence and the superharmonic minimum
principle give $w_j>M+1$ throughout these balls for large $j$.
The same lower bound for $W$ follows by comparison, using its blow-up
at the pole. All relevant minimizing points therefore lie in a fixed
regular compact set, where convergence is uniform. This proves
\eqref{eq:ub-minima-convergence}; in particular,
\begin{equation}\label{eq:ub-limit-normalization}
 \min_{\partial B_1}W=1.
\end{equation}

\subsection{Classification of the normalized harmonic limits}

Assume the bad sequence \eqref{eq:ub-bad-points} exists. By
\eqref{eq:ub-small-minimum}, at least one sequence of the preceding
type exists. Choose a centered sphere avoiding its concentration set.
The scaling law for the Pohozaev flux is
\[
 \mathcal P_1(u,r;0)=h(t_j)^p\mathcal P_{\eps_j}(w_j,R;0),
\]
where the left side is independent of its admissible radius $r$.
The rightmost flux is bounded by \eqref{eq:ub-harmonic-limit}.
It follows, before classifying the harmonic limit, that
\begin{equation}\label{eq:ub-original-zero-flux}
 \mathcal P_1(u,r;0)=0.
\end{equation}

Now consider any sequence with $h(t_j)\to0$ and any limit $W$
constructed above. The set $E=\Sigma\cup\{0\}$ is closed and
countable, hence has zero $p$-capacity. Since $W\ge0$, it extends to a
$p$-superharmonic function on $\R^n$.

If $\Sigma\ne\varnothing$, choose $z\in\Sigma$. Each $w_j$ solves
its critical equation in a full small ball about $z$, so
$\mathcal P_0(W,\rho;z)=0$. By Lemma~\ref{lem:ub-pole}, the constant
remainder at this pole is zero. For $0<\delta<k_z$, compare the
global superharmonic extension with
\[
 (k_z-\delta)(|x-z|^{-b}-R^{-b})
\]
on a large ball centered at $z$ with a small concentric ball removed.
The inner boundary inequality follows from the pole asymptotics and
the outer one from nonnegativity. Other exceptional points cause no
difficulty, because comparison is made with the superharmonic
extension. Shrinking the inner radius, letting $R\to\infty$, and
then letting $\delta\downarrow0$ gives
$W\ge k_z|x-z|^{-b}$. Lemma~\ref{lem:ub-pole} yields equality near
$z$. The contact set is open and closed in $\R^n\setminus E$:
at every contact point the comparison function has nonzero gradient,
so the local difference equation is uniformly elliptic.
The complement of a closed countable set in $\R^n$, $n\ge2$, is
connected. For example, two points can be joined through a third
point chosen off the countable union of lines determined by either
endpoint and an exceptional point. We conclude that
$W=k_z|x-z|^{-b}$ throughout the regular set. This excludes every
other pole in $\Sigma$.

If $\Sigma=\varnothing$, the only possible singularity of $W$ is
zero. If it is removable, the positive entire $p$-harmonic function
$W$ is constant: Harnack and the interior gradient estimate on
$B_{2R}(x)$ give $|\nabla W(x)|\le C W(x)/R$, and $R\to\infty$
forces its gradient to vanish. Otherwise \eqref{eq:ub-original-zero-flux} gives
the exact identity $\mathcal P_{\eps_j}(w_j,R;0)=0$ for every $j$ and fixed
admissible $R$. Passing to the limit and repeating the preceding
comparison argument gives $W=k_0|y|^{-b}$.

We have proved that, under the contradiction assumption, every
normalized limit associated with $h(t_j)\to0$ has one of the forms
\begin{equation}\label{eq:ub-limit-classification}
 W\equiv c>0,\qquad\text{or}\qquad W(y)=k|y-z|^{-b}
 \quad(k>0, z\in\R^n).
\end{equation}

\subsection{A second normalization and the contradiction}

Return to \eqref{eq:ub-bad-points}, and take a subsequence with
$x_j/r_j\to z\in\Sn$. In the normalization
\eqref{eq:ub-min-normalization} with $t_j=r_j$,
\[
 \Theta_j(x_j/r_j)=K_j^{1/a}\longrightarrow\infty.
\]
Hence $z\in\Sigma$. Equations \eqref{eq:ub-limit-classification}
and \eqref{eq:ub-limit-normalization} force
\[
 W(y)=2^b|y-z|^{-b}.
\]
By \eqref{eq:ub-minima-convergence},
\begin{equation}\label{eq:ub-first-minimum-profile}
 \frac{h(r_js)}{h(r_j)}\longrightarrow
 F(s):=\frac{2^bs^a}{(1+s)^b}
\end{equation}
uniformly on compact intervals of positive radii. The unique maximum
of $F$ occurs at $s=p-1$. Put $R_j=(p-1)r_j$. Then $h(R_j)\to0$
and
\begin{equation}\label{eq:ub-maximum-profile}
 \frac{h(R_js)}{h(R_j)}\longrightarrow
 f(s):=s^a\left(\frac p{1+(p-1)s}\right)^b.
\end{equation}
The function $f$ has a unique strict maximum at $s=1$.

Choose $1<c<2$ sufficiently close to one that
\[
 \min_{[1,c]}f>f(2),\qquad
 \min_{[1/c,1]}f>f(1/2).
\]
After taking a subsequence, assume $R_{j+1}<R_j/4$. Let $s_j$ be a
minimum point of $h$ on $[R_{j+1},R_j]$. The uniform convergence
in \eqref{eq:ub-maximum-profile}, applied at $R_{j+1}$ and compared
with $2R_{j+1}$, excludes $[R_{j+1},cR_{j+1}]$. Applied at $R_j$
and compared with $R_j/2$, it excludes $[R_j/c,R_j]$. Therefore
\[
 cR_{j+1}<s_j<R_j/c,\qquad h(s_j)\to0,\qquad
 h(s_jt)\ge h(s_j)\quad(c^{-1}\le t\le c).
\]
Normalize a second time, now at $s_j$. The compactness and
classification already proved apply, and
\eqref{eq:ub-minima-convergence} implies that
\[
 t^a\min_{|y|=t}W(y)\ge1\quad(c^{-1}\le t\le c),
 \qquad \min_{|y|=1}W(y)=1.
\]
Thus this function has a local minimum at $t=1$. But
\eqref{eq:ub-limit-classification} says it must be either $Ct^a$ or
\[
 \frac{kt^a}{(t+d)^b},\qquad d=|z|\ge0.
\]
The former is strictly increasing. The latter has derivative
\[
 kt^{a-1}(t+d)^{-b-1}\bigl(ad-(b-a)t\bigr).
\]
Since $b>a>0$, it is strictly decreasing when $d=0$, and when $d>0$
its only critical point is a strict maximum. Neither alternative
admits a local minimum. This contradiction proves that
$\sup_{0<|x|<r_0}|x|^au(x)<\infty$ for some $r_0>0$.
Continuity on the remaining compact annulus proves
Theorem~\ref{thm:upper}.

\section{Classification of isolated singularities}
\label{sec:classification-proof} 
In this section, we also use
\[
 \begin{gathered}
 p^*=\frac{np}{n-p},\qquad q=p^*-1,\qquad
 \sigma=q-p+1=\frac{p}{a},\\
 a=\frac{n-p}{p},\qquad b=\frac{n-p}{p-1}.
 \end{gathered}
\]
In particular, $b>a>0$.

\begin{theorem}[Removability criterion]\label{thm:removability}
Let $1<p<n$ and let $u\in W^{1,p}_{\mathrm{loc}}(B_1\setminus\{0\})$ be a nonnegative weak solution of \eqref{equ:main26}.
If
\begin{equation}\label{eq:rem-liminf}
 \liminf_{x\to0}|x|^{\frac{n-p}{p}} u(x)=0,
\end{equation}
then $u$ extends across the origin to a weak solution on $B_1$.
The extension belongs to $W^{1,p}_{\mathrm{loc}}(B_1)\cap C^{1,\gamma}_{\mathrm{loc}}(B_1)$ 
for some $\gamma\in(0,1)$.
\end{theorem}

We first record the classification of the harmonic profiles that will
arise after a normalization by spherical minima.  This is a
$p$-harmonic statement; see Ma-Zang \cite[Theorem 5.3]{MaZang}. 

\begin{lemma}[Positive punctured-space harmonic profiles]
\label{lem:rem-harmonic-profiles}
If $W>0$ is weakly $p$-harmonic in $\mathbb R^n\setminus\{0\}$, then
\begin{equation}\label{eq:rem-harmonic-classification}
 W(x)=L+\kappa|x|^{-b},\qquad L,\kappa\geq0.
\end{equation}
\end{lemma}

\begin{proof}[Proof of Theorem~\ref{thm:removability}]
If $u\equiv0$, there is nothing to prove.  Otherwise the strong
minimum principle gives $u>0$ in the punctured ball.  Local regularity
also gives $u\in C^{1,\gamma}_{\mathrm{loc}}(B_1\setminus\{0\})$.

\medskip\noindent
\emph{Normalization and annular compactness.}
For $0<r<R_0$, define
\begin{equation}\label{eq:rem-normalization}
 m(r)=\min_{|x|=r}u(x),\qquad
 h(r)=r^a m(r),\qquad
 w_r(y)=\frac{u(ry)}{m(r)}.
\end{equation}
The functions $m,h$ are positive and continuous.  A direct scaling
calculation gives
\begin{equation}\label{eq:rem-normalized-equation}
 -\Delta_p w_r=h(r)^\sigma w_r^q,
 \qquad \min_{|y|=1}w_r(y)=1.
\end{equation}
To obtain estimates independent of $m(r)$, write the original equation
as $-\Delta_pu=u^\sigma u^{p-1}$.  Its scaled potential satisfies
\[
 r^p u(ry)^\sigma\leq C_0^\sigma|y|^{-p}
\]
whenever $0<r|y|<R_0$, by $a\sigma=p$ and Theorem \ref{thm:upper}. Harnack inequalities for a bounded potential
on fixed annuli, followed by finite Harnack chains, therefore imply
\begin{equation}\label{eq:rem-spherical-harnack}
 \max_{|x|=r}u(x)\leq Hm(r)
\end{equation}
for all sufficiently small $r$, and, for every compact annulus
$K\Subset\mathbb R^n\setminus\{0\}$,
\begin{equation}\label{eq:rem-normalized-bounds}
 0<c_K\leq w_r\leq C_K\quad\text{on }K
\end{equation}
for sufficiently small $r$.  The constants are independent of $r$.
Since $h(r)\leq C_0$, the right-hand side in
\eqref{eq:rem-normalized-equation} is uniformly bounded on compact
annuli.  Local gradient regularity consequently gives uniform
$C^{1,\gamma}$ bounds there.  In particular, if
\begin{equation}\label{eq:rem-small-amplitudes}
 r_j\downarrow0,\qquad h(r_j)\longrightarrow0,
\end{equation}
then a subsequence satisfies
\begin{equation}\label{eq:rem-harmonic-limit}
 w_{r_j}\longrightarrow W
 \quad\text{in }C^1_{\mathrm{loc}}(\mathbb R^n\setminus\{0\}),
 \qquad -\Delta_pW=0,\quad W>0,\quad
 \min_{|y|=1}W=1.
\end{equation}
Positivity follows directly from the lower bounds in
\eqref{eq:rem-normalized-bounds}.

\medskip\noindent
\emph{The zero-flux restriction on harmonic limits.}
For $-\Delta_pv=\varepsilon v^{p^*-1}$, use the Pohozaev flux
\begin{equation}\label{eq:rem-pohozaev-flux}
 \begin{split}
 \mathcal P_\varepsilon(v,R)=\int_{\partial B_R}
 \bigg[&a v|\nabla v|^{p-2}v_\nu\\
 &+R\bigg(|\nabla v|^{p-2}v_\nu^2
       -\frac{|\nabla v|^p}{p}
       +\frac{\varepsilon v^{p^*}}{p^*}\bigg)\bigg]\,dS.
 \end{split}
\end{equation}
At $\nabla v=0$ the expressions are understood by continuity.
Lemma~\ref{lem:pohozaev} shows that this flux is independent of $R$
on concentric annuli; that lemma applies to the $C^1$ weak solutions
used here and does not require twice differentiability.  Denote the
constant flux of $u$ by $\mathcal P(u)$.  Scaling yields
\begin{equation}\label{eq:rem-flux-scaling}
 \mathcal P(u)=h(r)^p
       \mathcal P_{h(r)^\sigma}(w_r,1).
\end{equation}
Condition \eqref{eq:rem-liminf} gives a sequence as in
\eqref{eq:rem-small-amplitudes}.  The uniform $C^1$ bounds near
$\partial B_1$ imply that the normalized fluxes on the right-hand
side of \eqref{eq:rem-flux-scaling} remain bounded.  Hence
\begin{equation}\label{eq:rem-zero-flux}
 \mathcal P(u)=0.
\end{equation}
For every sequence satisfying \eqref{eq:rem-small-amplitudes},
passing to the limit in its normalized flux therefore gives
\[
 \mathcal P_0(W,1)=0.
\]

By Lemma~\ref{lem:rem-harmonic-profiles},
$W=L+\kappa|y|^{-b}$, where $L,\kappa\geq0$ and $L+\kappa=1$.
Substitution into \eqref{eq:rem-pohozaev-flux}, using
$(p-1)b/p=a$, gives
\[
 \mathcal P_0(W,1)
 =-a|\mathbb S^{n-1}|(b\kappa)^{p-1}L.
\]
It follows that every such normalized limit is one of precisely two
profiles:
\begin{equation}\label{eq:rem-two-profiles}
 W\equiv1\qquad\text{or}\qquad W(y)=|y|^{-b}.
\end{equation}

\medskip\noindent
\emph{From the lower limit to a full limit.}
We claim that
\begin{equation}\label{eq:rem-amplitude-limit}
 h(r)\longrightarrow0\qquad \text{as } r\downarrow0.
\end{equation}
Suppose instead that $\limsup_{r\downarrow0}h(r)>0$ and choose
$d>0$ smaller than half this upper limit.  By continuity and
$\liminf_{r\downarrow0}h(r)=0$, there are intervals
$[\alpha_j,\beta_j]$ with $\beta_j\downarrow0$, whose endpoint
values are at least $2d$, and with an interior point where
$h<1/j$.  Indeed, choose successively a sufficiently small outer
radius with $h\geq2d$, a smaller radius with $h<1/j$, and a still
smaller radius with $h\geq2d$.  For all large $j$, a minimum point
$s_j$ of $h$ on this interval lies in its interior and satisfies
$h(s_j)\to0$.  Thus $s_j$ is a local minimum of $h$.

Passing to a subsequence, $w_{s_j}$ converges in $C^1$ on compact
annuli to one of the two functions in \eqref{eq:rem-two-profiles}.
Set
\[
 G_j(t,\theta)=t^a w_{s_j}(t\theta),
 \qquad t>0,\quad\theta\in\mathbb S^{n-1}.
\]
Then
\[
 \partial_tG_j(t,\theta)
 =t^{a-1}\bigl(a w_{s_j}(t\theta)
              +t\theta\cdot\nabla w_{s_j}(t\theta)\bigr).
\]
In a fixed neighborhood of $t=1$, these derivatives are uniformly
strictly positive for all large $j$ if the limit is $1$, and
uniformly strictly negative if the limit is $|y|^{-b}$.  The
corresponding limiting derivatives are
$a t^{a-1}>0$ and $(a-b)t^{a-b-1}<0$, respectively.  Taking minima
over $\theta$ preserves this strict monotonicity, since the derivative
has a uniform strict sign.  But
\[
 \min_\theta G_j(t,\theta)=\frac{h(s_jt)}{h(s_j)}
\]
has a local minimum at $t=1$, a contradiction.  Notice that no
differentiability of the spherical minimum $m$ has been used.
This proves \eqref{eq:rem-amplitude-limit}; by
\eqref{eq:rem-spherical-harnack}, it also gives
\begin{equation}\label{eq:rem-uniform-critical-decay}
 \sup_{|x|=r}r^a u(x)\longrightarrow0\qquad \text{as } r\downarrow0.
\end{equation}

\medskip\noindent
\emph{Elimination of the fundamental-solution profile.}
Define
\[
 F(r,\theta)=\frac{a u(r\theta)+r\,\partial_ru(r\theta)}{m(r)},
 \qquad c=\frac12\min\{a,b-a\}>0.
\]
The compactness already established, together with
\eqref{eq:rem-amplitude-limit} and \eqref{eq:rem-two-profiles}, implies
that there exists $R>0$ such that every $r\in(0,R)$ satisfies exactly
one of the alternatives
\begin{equation}\label{eq:rem-sign-alternatives}
 \begin{aligned}
 &F(r,\theta)\geq c\quad\text{for all }\theta\in\mathbb S^{n-1},
 \qquad\text{or}\\
 &F(r,\theta)\leq-c\quad\text{for all }\theta\in\mathbb S^{n-1}.
 \end{aligned}
\end{equation}
Otherwise a sequence of exceptional radii tending to zero would
have a normalized subsequence converging to one of the two profiles;
on the unit sphere the corresponding functions $F$ converge
uniformly to $a$ or $a-b$, contradicting exceptionality.

Fix $\theta_0\in\mathbb S^{n-1}$.  The function
$r\mapsto F(r,\theta_0)$ is continuous and, by
\eqref{eq:rem-sign-alternatives}, stays away from zero.  Its sign is
therefore constant on $(0,R)$.  The negative sign is impossible:
it would imply
\[
 \frac{d}{dr}\bigl(r^a u(r\theta_0)\bigr)
 =r^{a-1}m(r)F(r,\theta_0)<0,
\]
and hence $r^a u(r\theta_0)>r_0^a u(r_0\theta_0)>0$ whenever
$0<r<r_0<R$, contrary to \eqref{eq:rem-uniform-critical-decay}.
Thus only the positive alternative occurs.  A subsequential
fundamental-solution limit would make $F$ converge to $a-b<0$,
so that profile is excluded.  Compactness now gives convergence of
the entire family:
\begin{equation}\label{eq:rem-constant-normalized-limit}
 w_r\longrightarrow1
 \quad\text{in }C^1_{\mathrm{loc}}(\mathbb R^n\setminus\{0\})
 \quad\text{as }r\downarrow0.
\end{equation}
In particular,
\begin{equation}\label{eq:rem-small-log-gradient}
 \sup_{|x|=r}\frac{r|\nabla u(x)|}{u(x)}
 \leq\sup_{|x|=r}\frac{r|\nabla u(x)|}{m(r)}
 \longrightarrow0.
\end{equation}

\medskip\noindent
\emph{Arbitrarily small power growth and finite energy.}
Given $\delta>0$, \eqref{eq:rem-small-log-gradient} implies
$|\partial_r\log u(r\theta)|\leq\delta/r$ for all sufficiently
small $r$, uniformly in $\theta$.  Integrating along rays up to a
fixed sufficiently small sphere and using its finite maximum gives
$u(x)\leq C_\delta|x|^{-\delta}$.  Applying
\eqref{eq:rem-small-log-gradient} once more yields, after adjusting
the constant,
\begin{equation}\label{eq:rem-arbitrary-small-growth}
 u(x)\leq C_\delta|x|^{-\delta},\qquad
 |\nabla u(x)|\leq C_\delta|x|^{-1-\delta}, 
 \qquad 0<|x|<R_\delta.
\end{equation}
Choose
\begin{equation}\label{eq:rem-delta-choice}
 0<\delta<\min\{a,p/q\},
\end{equation}
and fix $R<R_\delta$.  Since $p(1+\delta)<n$,
\eqref{eq:rem-arbitrary-small-growth} implies that $u$ and its
punctured-ball gradient belong to $L^p(B_R)$.

To verify the Sobolev extension without assuming it, let
$\chi_\varepsilon$ be a smooth cutoff that is zero on $B_\varepsilon$,
one outside $B_{2\varepsilon}$, takes values in $[0,1]$, and satisfies
$|\nabla\chi_\varepsilon|\leq C/\varepsilon$.
The function $\chi_\varepsilon u$, extended by zero near the origin,
belongs to $W^{1,p}(B_R)$.  Moreover,
\[
 \int_{B_{2\varepsilon}\setminus B_\varepsilon}
       |u\nabla\chi_\varepsilon|^p
 \leq C\varepsilon^{n-p-p\delta}\longrightarrow0.
\]
Consequently $\chi_\varepsilon u\to u$ in $L^p(B_R)$, and its
gradients converge in $L^p(B_R)$ to the original gradient, extended
arbitrarily at the origin.  Closedness of the weak derivative proves
\begin{equation}\label{eq:rem-sobolev-extension}
 u\in W^{1,p}(B_R).
\end{equation}

\medskip\noindent
\emph{Extension of the equation.}
For $\varphi\in C_c^\infty(B_R)$, test the punctured-ball equation
with $\varphi\chi_\varepsilon$.  The additional term is bounded by
\[
 \begin{split}
 \left|\int_{B_R}\varphi|\nabla u|^{p-2}\nabla u
                       \cdot\nabla\chi_\varepsilon\right|
 &\leq C\|\varphi\|_\infty
 \left(\int_{B_{2\varepsilon}\setminus B_\varepsilon}
                       |\nabla u|^p\right)^{(p-1)/p}
 \varepsilon^{(n-p)/p}\longrightarrow0.
 \end{split}
\]
Also, $u^q\in L^1(B_R)$ because $\delta q<p<n$.  Dominated
convergence in the remaining terms therefore yields
\begin{equation}\label{eq:rem-extended-equation}
 \int_{B_R}|\nabla u|^{p-2}\nabla u\cdot\nabla\varphi
   =\int_{B_R}u^q\varphi
 \qquad(\varphi\in C_c^\infty(B_R)).
\end{equation}
Thus no point source is created at the origin.

\medskip\noindent
\emph{Boundedness and regularity at the origin.}
It remains to improve the Sobolev extension to a continuous one.
Put
\[
 s=\delta q<p,\qquad \alpha=\frac{p-s}{p-1}>0.
\]
By \eqref{eq:rem-arbitrary-small-growth},
$u^q(x)\leq C|x|^{-s}$ in $B_R\setminus\{0\}$.  Choose
$M\geq\max_{\partial B_R}u$ and define
\[
 Z(x)=M+K(R^\alpha-|x|^\alpha).
\]
For $x\neq0$, direct calculation gives
\begin{equation}\label{eq:rem-bounded-barrier}
 -\Delta_pZ=(K\alpha)^{p-1}(n-s)|x|^{-s}.
\end{equation}
Since $\alpha>0$ and $n>p$, the gradient of $Z$ belongs to $L^p(B_R)$,
so $Z\in W^{1,p}(B_R)\cap C(\overline{B_R})$.  Its flux through
$\partial B_\varepsilon$ is $O(\varepsilon^{n-s})$, which tends to
zero.  Integrating by parts outside $B_\varepsilon$ and passing to
the limit therefore proves \eqref{eq:rem-bounded-barrier} weakly on
the whole ball.  Choose $K$ sufficiently large that its right-hand
side dominates $u^q$ almost everywhere.

The comparison test is legitimate in the whole ball.  Indeed, if
$t=(p^*)'=np/[n(p-1)+p]$, then $pt<n$ because $p<n$, and therefore
$st<n$.  Both sources belong to $L^{(p^*)'}(B_R)$; the weak equations
extend by density to tests in $W^{1,p}_0(B_R)$.  The boundary trace
inequality $u\leq Z$ on $\partial B_R$ gives
$(u-Z)_+\in W^{1,p}_0(B_R)$.  Subtracting the equations and using this
test yields
\[
 \int_{\{u>Z\}}
 \bigl(|\nabla u|^{p-2}\nabla u-|\nabla Z|^{p-2}\nabla Z\bigr)
       \cdot(\nabla u-\nabla Z)\leq0.
\]
Strict monotonicity of $\xi\mapsto|\xi|^{p-2}\xi$, followed by
the zero boundary trace, implies $(u-Z)_+=0$.  Thus
\[
 u\leq Z\leq M+KR^\alpha\qquad\text{almost everywhere in }B_R.
\]
The right-hand side $u^q$ of \eqref{eq:rem-extended-equation} is now
bounded.  Standard interior gradient regularity for the
$p$-Laplacian gives a representative in
$C^{1,\gamma}_{\mathrm{loc}}(B_R)$, for some $\gamma\in(0,1)$.
This representative agrees with the original continuous function
off the origin. This proves the theorem.  
\end{proof}

\begin{proof}[Proof of Theorem~\ref{thm:main}]
If $u\equiv0$, alternative~\textup{(i)} holds. Otherwise the strong minimum principle gives $u>0$. Theorem~\ref{thm:upper} gives the upper bound on $B_{1/2}\setminus\{0\}$. If $\liminf_{x\to0}|x|^{\frac{n-p}{p}} u(x)=0$, Theorem~\ref{thm:removability} gives alternative~\textup{(i)}.
If the origin is nonremovable, the contrapositive of Theorem~\ref{thm:removability} gives
\[
 \liminf_{x\to0}|x|^{\frac{n-p}{p}} u(x)>0.
\]
Together with the upper bound this yields the two-sided estimate near zero. Positivity and continuity on the remaining compact annulus extend it to $0<|x|<1/2$. 
This proves Theorem~\ref{thm:main}. 
\end{proof}

\section{Radial rigidity under two-sided critical bounds}
\label{sec:global-symmetry}

Throughout this section,
\[
 a=\frac{n-p}{p},\qquad p^*=\frac{np}{n-p},\qquad q=p^*-1,\qquad1<p<n.
\]
Constants in the two-sided bounds may depend on the solution.
We use continuous Steiner symmetrization in the sense of
Brock~\cite{Brock2000}, with the localization argument specified below.

\subsection{Quantitative continuous Steiner dissipation}

Fix a direction and write $x=(z,y)\in\mathbb R^{n-1}\times\mathbb R$.
For $0\leq w\in C_c^1(\mathbb R^n)$, write $w^\tau$ for its continuous
Steiner symmetrization.  Until first contact with another component,
an interval component $(l,r)$ of a fiber superlevel set evolves by
\begin{equation}\label{eq:cst-interval-motion}
 l^\tau=e^{-\tau}c-h,\qquad r^\tau=e^{-\tau}c+h,\qquad
 c=\frac{l+r}{2},\quad h=\frac{r-l}{2}.
\end{equation}
Intervals merge at contact and thereafter obey the same rule.
The transformation is equimeasurable, nonexpansive in $L^p$, preserves
Lipschitz bounds, and decreases the Dirichlet integral; see
\cite[Remark~2.6(2), (5), (7), and Theorem~3.2]{Brock2000}.

An interval component $(l(z,s),r(z,s))$ of $\{y:w(z,y)>s\}$ is a
\emph{regular pair} if $w_y(z,l)>0$ and $w_y(z,r)<0$.
Its endpoints are $C^1$ functions of $(z,s)$ in a regular parameter
neighborhood.  Integration over the collection $\mathcal P_w$ of
regular pairs means summation over interval components followed by
integration in $dz\,ds$, counting each pair once.  Equivalently,
partition the collection into disjoint measurable pieces contained
in regular parameter charts.  Put
\begin{align*}
 c&=(l+r)/2,&h&=(r-l)/2,\\
 A&=(-\nabla_zh,-h_s),&B&=(\nabla_zc,c_s),\\
 X_l&=(\nabla_zl,l_s)=A+B,&X_r&=(-\nabla_zr,-r_s)=A-B,
\end{align*}
and define, for $\lambda>0$,
\begin{equation}\label{eq:cst-F-definition}
 F(\zeta,\lambda)=(1+|\zeta|^2)^{p/2}\lambda^{1-p},\qquad
 \mathcal D_w=(DF(X_l)-DF(X_r))\cdot B.
\end{equation}
The last components of $X_l,X_r$ are positive.  The function $F$
is strictly convex: it is the perspective
$\lambda|(\zeta,1)/\lambda|^p$ of $|\cdot|^p$, and equality forces
both arguments to agree because the last component of $(\zeta,1)$
is fixed.  Hence
\begin{equation}\label{eq:cst-defect-zero}
 \mathcal D_w\geq0,\qquad\mathcal D_w=0\ \Longleftrightarrow\ B=0.
\end{equation}

The following quantitative refinement of Brock's local energy
argument is proved here.

\begin{lemma}[Quantitative localized dissipation]
\label{lem:cst-quantitative}
Let $0\leq w\in C_c^1(\mathbb R^n)$ and set
\[
 \mathcal E(v)=\frac1p\int_{\mathbb R^n}|\nabla v|^p\,dx.
\]
Then
\begin{equation}\label{eq:cst-global-defect}
 \liminf_{\tau\downarrow0}
 \frac{\mathcal E(w)-\mathcal E(w^\tau)}{\tau}
 \geq\frac1p\int_{\mathcal P_w}\mathcal D_w\,dz\,ds.
\end{equation}
The right side is initially a nonnegative extended integral.
Moreover, $\dot w=\lim_{\tau\downarrow0}(w^\tau-w)/\tau$ exists
almost everywhere, equals $cw_y$ at regular endpoints, and is zero
wherever $w_y=0$.
\end{lemma}

\begin{proof}
Choose finitely many relatively compact regular parameter patches
$Q_1,\ldots,Q_N$ whose endpoint images have pairwise disjoint closures.
Let $G_j(\tau)$ be the union of their two endpoint images under
$(z,s)\mapsto(z,l^\tau(z,s))$ and $(z,s)\mapsto(z,r^\tau(z,s))$,
and set $G(\tau)=\bigcup_jG_j(\tau)$.
Uniform endpoint regularity gives a positive empty gap immediately
outside each selected interval.  Finite propagation gives a common
$\tau_0>0$ such that none of these intervals merges for
$0\leq\tau<\tau_0$ and their endpoint images remain disjoint.
For these times,
\begin{equation}\label{eq:cst-finite-complement}
 \int_{\mathbb R^n\setminus G(0)}|\nabla w|^p\,dx
 \geq\int_{\mathbb R^n\setminus G(\tau)}|\nabla w^\tau|^p\,dx.
\end{equation}
This is the finite-patch version of (6.13) in~\cite{Brock2000}.
Here is why that extension holds.  For the \emph{good} piecewise
smooth functions used there, the proof of \cite[Theorem~3.1]{Brock2000} decomposes
the coarea energy into interval contributions.  Its inequality
(3.16) applies to each unmerged interval separately, and (3.17)
applies to each merger block separately.  At every $(z,s)$ delete
the selected surviving pairs from the sum.  All remaining merger
blocks consist entirely of unselected intervals, so the same
inequalities give \eqref{eq:cst-finite-complement}.
Thus no overlapping complement inequalities are added.

For completeness, the required relative approximation can be constructed
as follows. Let $K$ be the union of the closures of the selected endpoint
images. Choose $\delta>0$ and finitely many open boxes $V_\nu$ with
$\overline{V_\nu}\subset\{|w_y|>2\delta\}$ whose smaller concentric
boxes cover $K$. Choose smooth functions $0\leq\eta_\nu\leq1$ which
equal one on those smaller boxes and are positive precisely in $V_\nu$,
and put $\chi=1-\prod_\nu(1-\eta_\nu)$. Thus $\chi=1$ on a
neighborhood of $K$, while $\{\chi>0\}$ is a finite union of boxes
contained in $\{|w_y|>2\delta\}$.

Nonnegative nodal piecewise affine interpolation on successively finer
shape-regular simplicial triangulations gives compactly supported functions $b_j$ with
common bounded support such that
\[
 \|b_j-w\|_\infty+
 \mathop{\rm ess\,sup}|\nabla b_j-\nabla w|\longrightarrow0.
\]
Fix $R$ large enough that these supports are compactly contained in the
positivity set of the piecewise affine tent
\[
 \rho(z,y)=\left(R-|y|-\sum_{i=1}^{n-1}|z_i|\right)_+.
\]
On the positive affine pieces of $\rho$, its vertical derivative is
$1$ or $-1$. For each $j$, the common finite refinement of $b_j$ and
$\rho$ has only finitely many vertical slopes. We may therefore choose
$\varepsilon_j\downarrow0$ avoiding the finitely many values for which
one of the positive affine pieces of
$a_j=b_j+\varepsilon_j\rho$ has zero vertical slope. The functions
$a_j$ are nonnegative, have common bounded support, and converge to $w$
with the same uniform function and essential gradient convergence.

Set $w_j=\chi w+(1-\chi)a_j$. These functions agree with $w$ near $K$,
remain nonnegative, and converge uniformly and strongly in $W^{1,p}$.
On $\{\chi>0\}$,
\[
 (w_j)_y=w_y+(1-\chi)((a_j)_y-w_y)+\chi_y(w-a_j).
\]
For sufficiently large $j$ this derivative has the sign of $w_y$ and
absolute value at least $\delta$. Each vertical fiber meets
$\{\chi>0\}$ in finitely many intervals; on each such interval $w_y$
has a constant sign. Outside this set, $w_j=a_j$ is piecewise affine
with finitely many pieces and nonzero slope wherever positive. Hence
$w_j$ has finitely many intersections with each positive fiber level,
and its directional slopes on positive pieces have a positive lower
bound in absolute value for each $j$. It is therefore a good piecewise
$C^1$ function in the sense needed
for Brock's inverse-coordinate argument.

The selected endpoint graphs agree exactly with those of $w$.
Between their endpoint neighborhoods, compactness and $w>s$ give
a uniform positive level gap.  Uniform convergence therefore
preserves the selected interval components; the exterior gaps give
a common nonmerging time.  Their moving sets $G(\tau)$ are consequently
independent of $j$.  Apply \eqref{eq:cst-finite-complement} to $w_j$.
The left side converges strongly.  Nonexpansivity gives
$w_j^\tau\to w^\tau$ in $L^p$, and the Dirichlet inequality gives
weak $W^{1,p}$ compactness.  Lower semicontinuity on the fixed set
$\mathbb R^n\setminus G(\tau)$ proves the inequality for $w$.
All directionally critical fibers remain in the complement, without
any assumption that their energy vanishes.

On a selected patch, change of variables on both endpoint branches gives
\[
 \int_{G_j(\tau)}|\nabla w^\tau|^p\,dx
 =\int_{Q_j}[F(A+e^{-\tau}B)+F(A-e^{-\tau}B)]\,dz\,ds.
\]
Differentiation at zero is valid on this compact regular patch and
gives $-\int_{Q_j}\mathcal D_w$.
Combine this with \eqref{eq:cst-finite-complement}, then exhaust the
regular-pair collection by finite disjoint parameter pieces compactly
contained in regular charts.  These may be chosen on the regular
left-endpoint set; injectivity of the pairing map makes their right
images disjoint as well.  First shrinking finitely many pieces gives
disjoint closures, and inner approximation recovers omitted boundaries.
Nonnegativity of the defect proves \eqref{eq:cst-global-defect}.

For the first variation, suppose $\operatorname{supp}w\subset B_R$.
On each fiber every interval remains in $[-R,R]$, and its endpoints
move with speed at most $R$. At a merger the exterior endpoints are
continuous, and the merged interval again has speed at most $R$.
Mergers remove interior boundaries and create no new gaps.
For a finite union of intervals $E$, it follows that
\[
 E\mathbin{\ominus}[-R\tau,R\tau]\subset E^\tau
 \subset E+[-R\tau,R\tau],
\]
where $E\mathbin{\ominus}[-d,d]=\{y:[y-d,y+d]\subset E\}$.
For an arbitrary bounded open fiber set, approximate from inside
by finite unions of intervals. Monotonicity and continuity from
inside pass these inclusions to the limit. Applying them to the
strict superlevel sets of $w$ and then taking limits in the level
gives, for $0<\tau<1$,
\begin{equation}\label{eq:cst-finite-propagation}
 \inf_{|t|\leq R\tau}w(z,y+t)\leq w^\tau(z,y)
 \leq\sup_{|t|\leq R\tau}w(z,y+t).
\end{equation}
Thus $\|w^\tau-w\|_\infty\leq C_w\tau$, with a common compact
support.  If $w_y=0$, its $C^1$ expansion along the fiber gives
$\dot w=0$.  At regular endpoints, \eqref{eq:cst-interval-motion}
and the inverse identity give $\dot w=cw_y$.
For each fixed $z$, the critical values of the $C^1$ function
$y\mapsto w(z,y)$ form a null set.  Their inverse images in
$\{w_y\ne0\}$ have measure zero by local change of variables.
Fubini shows that almost every point with $w_y\ne0$ belongs to a
regular pair, proving the remaining claims.
\end{proof}

\subsection{Rigidity under the critical two-sided bound}

\begin{proposition}[Radial rigidity under two-sided critical bounds]
\label{prop:two-sided-rigidity}
Let $U\in W^{1,p}_{\mathrm{loc}}(\mathbb R^n\setminus\{0\})$ be a positive weak solution of
\begin{equation}\label{eq:global-critical-equation}
 -\Delta_pU=U^{p^*-1} \quad\text{in }\mathbb R^n\setminus\{0\}.
\end{equation}
Suppose $0<c_0\leq C_0<\infty$ and
\begin{equation}\label{eq:global-critical-two-sided}
 c_0|x|^{-\frac{n-p}{p}}\leq U(x)\leq C_0|x|^{-\frac{n-p}{p}}, \qquad x\ne0.
\end{equation}
Then $U(x)=\phi(|x|)$ for some $\phi\in C^1((0,\infty))$ with
$\phi'(r)<0$ for every $r>0$.
\end{proposition}

\begin{proof}
Scaling and interior $C^{1,\alpha}$ estimates give
\begin{equation}\label{eq:global-critical-gradient}
 U\in C^{1,\alpha}_{\mathrm{loc}}(\mathbb R^n\setminus\{0\}),
 \qquad|\nabla U(x)|\leq C|x|^{-a-1}.
\end{equation}
Indeed, $r^aU(r\,\cdot)$ solves the same equation and is uniformly
bounded on fixed compact annuli for every $r>0$.
Since $p(a+1)=n$, for $T\geq1$,
\begin{equation}\label{eq:global-log-energy}
 \begin{split}
 \int_{\{e^{-2T}<U<e^{2T}\}}|\nabla U|^p\,dx
 &\leq C\int_{(c_0e^{-2T})^{1/a}<|x|<
                    (C_0e^{2T})^{1/a}}|x|^{-n}\,dx\\
 &\leq CT.
 \end{split}
\end{equation}
Constants denoted by $C$ below are independent of $T$ and the direction.

Fix any direction $y$.  For defining superlevel sets only, set
$U(0)=+\infty$.  Every positive superlevel set is bounded.
On its regular fiber components use the preceding notation
$l,r,c,h,A,B,X_l,X_r$, with $s$ the level of $U$.
Put
\[
 \mathcal D=(DF(X_l)-DF(X_r))\cdot B,\qquad
 \lambda_l=l_s,\quad\lambda_r=-r_s,\quad
 g_l=|\nabla U(z,l)|,\quad g_r=|\nabla U(z,r)|,
\]
and
\[
 \mathcal Q=F(X_l)+F(X_r),\qquad
 \lambda_*=\min\{\lambda_l,\lambda_r\},\qquad\mathcal P=\mathcal P_U.
\]
Implicit differentiation and one-dimensional change of variables give
\begin{equation}\label{eq:global-pair-energy}
 \begin{gathered}
 F(X_l)=\lambda_lg_l^p,\qquad F(X_r)=\lambda_rg_r^p,\\
 \int_{\mathcal P\cap\{\alpha<s<\beta\}}\mathcal Q\,dz\,ds
 \leq\int_{\{\alpha<U<\beta\}}|\nabla U|^p\,dx
 \qquad(0<\alpha<\beta<\infty).
 \end{gathered}
\end{equation}
Only an inequality is needed; directionally critical sets may carry
additional energy.

For $g=(1+|\zeta|^2)^{1/2}/\lambda$,
\[
 F_\lambda=-(p-1)g^p,\qquad F_{\lambda\lambda}=p(p-1)F/\lambda^2.
\]
Hessian Cauchy--Schwarz along $X(t)=A+tB$, $-1\leq t\leq1$, gives
\begin{align*}
 |F_\lambda(X_l)-F_\lambda(X_r)|^2
 &\leq
 \left(\int_{-1}^1B^TD^2F(X(t))B\,dt\right)
 \left(\int_{-1}^1F_{\lambda\lambda}(X(t))\,dt\right)\\
 &\leq\frac{p(p-1)}{\lambda_*^2}\mathcal D\mathcal Q.
\end{align*}
The first factor equals $\mathcal D$, while convexity implies
$\int_{-1}^1F(A+tB)\,dt\leq F(X_l)+F(X_r)$.
Consequently,
\begin{equation}\label{eq:global-flux-bregman}
 (g_l^p-g_r^p)^2
 \leq\frac p{p-1}\frac{\mathcal D\mathcal Q}{\lambda_*^2}.
\end{equation}
On $U=s$, the two-sided and gradient bounds give
\[
 |c|\leq Cs^{-1/a},\qquad g_l+g_r\leq Cs^{1+1/a},
 \qquad\lambda_l,\lambda_r\geq C^{-1}s^{-1-1/a}.
\]
Thus $|c|/s\leq C\lambda_*$.  With
$\mathcal J=(c/s)(g_l^p-g_r^p)$, we obtain
\begin{equation}\label{eq:global-flux-defect}
 |\mathcal J|^2\leq C\mathcal D\mathcal Q.
\end{equation}

Choose $\chi\in C_c^\infty(\mathbb R)$ with
\[
 0\leq\chi\leq1,\qquad \chi=1\text{ on }[-1,1],\qquad
 \chi>0\text{ on }(-2,2),\qquad\operatorname{supp}\chi=[-2,2].
\]
Fix an integer $m\geq\max\{2,2/(p-1)\}$.  For $T\geq1$ define
\begin{equation}\label{eq:global-value-cutoff}
 \begin{aligned}
 \Psi_T'(s)&=\chi(\log s/T)^m,
 &\Psi_T(s)&=\int_0^s\Psi_T'(t)\,dt,\\
 \eta_T(\ell)&=\chi(\ell/T)^{mp},
 &\Psi_T'(s)^p&=\eta_T(\log s).
 \end{aligned}
\end{equation}
The derivative is understood to be zero at $s=0$.
Let $w=\Psi_T(U)$, extended by its constant limiting value at the
origin.  It is constant near zero, vanishes outside a sufficiently
large ball, and belongs to $C_c^1(\mathbb R^n)$.
The weak chain rule gives on all of $\mathbb R^n$
\begin{equation}\label{eq:global-truncated-equation}
 -\Delta_pw
 =\Psi_T'(U)^{p-1}U^q
 -(p-1)\Psi_T'(U)^{p-2}\Psi_T''(U)|\nabla U|^p
 =:f_T.
\end{equation}
To verify this identity, test the equation for $U$ with
$\Psi_T'(U)^{p-1}\varphi$, which is supported away from the origin
and infinity.  The expression $\Psi_T'^{\,p-2}\Psi_T''$ extends
continuously by zero across the cutoff endpoints: its factor involving
$\chi$ is a constant multiple of
$\chi^{m(p-1)-1}\chi'$, and $m(p-1)-1\geq1$.
Thus $f_T$ is bounded and supported in a compact annulus, including
when $1<p<2$.

Apply Lemma~\ref{lem:cst-quantitative}.  On the active level range,
write $\varsigma=\Psi_T(s)$.  The inverse graphs for $w$ and $U$ have
the same $z$ derivatives, whereas their level derivatives differ by
the factor $1/\Psi_T'(s)$.  Homogeneity of $F$ in its last variable
therefore makes the defect for $w$ equal to
$\Psi_T'(s)^{p-1}\mathcal D$.  Multiplication by
$d\varsigma=\Psi_T'(s)\,ds$ gives
\begin{equation}\label{eq:global-truncated-dissipation}
 \liminf_{\tau\downarrow0}
 \frac{\mathcal E(w)-\mathcal E(w^\tau)}{\tau}
 \geq\frac1p\int_{\mathcal P}\eta_T(\log s)\mathcal D\,dz\,ds.
\end{equation} 
Convexity of the energy and testing
\eqref{eq:global-truncated-equation} with $w^\tau-w$ yield
\begin{equation}\label{eq:global-energy-convexity}
 \mathcal E(w^\tau)-\mathcal E(w)
 \geq\int_{\mathbb R^n}f_T(w^\tau-w)\,dx.
\end{equation}
This is an admissible compactly supported $W^{1,p}$ test, by the
Lipschitz and support properties of continuous symmetrization and
density of smooth tests.
Since $\|w^\tau-w\|_\infty\leq C_T\tau$ and $f_T\in L^\infty_c$,
\begin{equation}\label{eq:global-finite-loss}
 \limsup_{\tau\downarrow0}
 \frac{\mathcal E(w)-\mathcal E(w^\tau)}{\tau}<\infty.
\end{equation}
Together with the finite-patch proof of
\eqref{eq:global-truncated-dissipation}, this already establishes
\begin{equation}\label{eq:global-AT-definition}
 A_T:=\int_{\mathcal P}\eta_T(\log s)\mathcal D\,dz\,ds<\infty.
\end{equation}
Finiteness is thus known before any weighted absorption is performed.

The first term of $f_T$ has zero first variation.
On the strictly increasing range of $\Psi_T$ define
\[
 K_T(\Psi_T(s))=\Psi_T'(s)^{p-1}s^q.
\]
It extends continuously by zero at the two flat endpoint values.
For $H_T(t)=\int_0^tK_T(v)\,dv$, equimeasurability gives
$\int H_T(w^\tau)=\int H_T(w)$.
Uniform continuity of $K_T$ on this bounded range, the bound
$\|w^\tau-w\|_\infty\leq C_T\tau$, and the common compact support
imply
\begin{equation}\label{eq:global-source-cancellation}
 \lim_{\tau\downarrow0}\frac1\tau
 \int_{\mathbb R^n}K_T(w)(w^\tau-w)\,dx=0.
\end{equation}

For the second term, dominated convergence and
Lemma~\ref{lem:cst-quantitative} give
\begin{align}
 \limsup_{\tau\downarrow0}
 \frac{\mathcal E(w)-\mathcal E(w^\tau)}{\tau}
 &\leq(p-1)\int_{\mathbb R^n}
 \Psi_T'(U)^{p-1}\Psi_T''(U)c\,U_y|\nabla U|^p\,dx
 \notag\\
 &=\frac{p-1}{p}\int_{\mathcal P}
       \eta_T'(\log s)\mathcal J\,dz\,ds.
 \label{eq:global-chain-variation}
\end{align}
A prime on $\eta_T$ denotes differentiation in its argument $\ell$.
The equality uses
\[
 \Psi_T'(s)^{p-1}\Psi_T''(s)=\frac{\eta_T'(\log s)}{ps}
\]
and change of variables on the increasing and decreasing fiber
branches.  At directionally critical points $U_y=0$; exceptional
critical levels are null as in the lemma.  No contribution is
omitted.  For fixed $T$ the integral is absolutely convergent:
the spatial integrand is bounded and supported in a compact annulus. 
Combining \eqref{eq:global-truncated-dissipation} and
\eqref{eq:global-chain-variation} yields
\[
 A_T\leq(p-1)\int_{\mathcal P}\eta_T'(\log s)\mathcal J\,dz\,ds.
\]
By \eqref{eq:global-flux-defect} and Young's inequality,
\[
 A_T\leq C\int_{\mathcal P}|\eta_T'|\sqrt{\mathcal D\mathcal Q}
 \leq\frac12A_T+
 C\int_{\mathcal P}\frac{|\eta_T'|^2}{\eta_T}\mathcal Q\,dz\,ds.
\]
Define the quotient as zero where $\eta_T=0$.  Since $mp\geq2$,
\[
 \frac{|\eta_T'(\ell)|^2}{\eta_T(\ell)}
 =\frac{m^2p^2}{T^2}\chi(\ell/T)^{mp-2}|\chi'(\ell/T)|^2
 \leq\frac C{T^2},
\]
and it vanishes outside $[-2T,2T]$.
Equations \eqref{eq:global-pair-energy} and
\eqref{eq:global-log-energy} prove
\begin{equation}\label{eq:global-vanishing-defect}
 \int_{\mathcal P}\eta_T(\log s)\mathcal D\,dz\,ds
 \leq\frac C{T^2}\int_{\{e^{-2T}<U<e^{2T}\}}|\nabla U|^p\,dx
 \leq\frac CT.
\end{equation}

On every fixed relatively compact regular parameter patch,
$\eta_T(\log s)=1$ for all sufficiently large $T$.
Letting $T\to\infty$ proves that $\mathcal D$ has zero integral on
the patch.  Its continuity follows from the $C^1$ endpoint graphs.
Strict convexity of $F$ therefore gives, pointwise,
\[
 \mathcal D=0,\qquad c_s=0,\qquad\nabla_zc=0.
\]
Implicit differentiation translates this into
\begin{equation}\label{eq:global-paired-derivatives}
 U_y(z,l)=-U_y(z,r),\qquad
 \nabla_zU(z,l)=\nabla_zU(z,r).
\end{equation}

These identities also exclude a critical opposite endpoint.
Suppose $U_y(z,l)>0$ and $(l,r)$ is the component of
$\{U(z,\cdot)>s\}$ adjoining $l$.
Choose regular values $s_j\downarrow s$ of the one-dimensional
function.  Its possible singular point at zero does not affect
finite levels: one-dimensional Sard applies on either side.
The corresponding component endpoints $l_j,r_j$ tend to $l,r$.
For the left endpoint this follows from the implicit function
theorem, and for the right endpoint from the strict inequality
$U>s$ on each compact subinterval of $(l,r)$.
For large $j$, $U_y(z,l_j)>0$ and $U_y(z,r_j)<0$, so
\eqref{eq:global-paired-derivatives} applies.
Continuity gives $U_y(z,r)=-U_y(z,l)<0$ and the transverse identity.
Thus the paired-endpoint identities hold at every finite positive
level of $U$, in every direction, whenever the entering endpoint
has positive directional derivative.

We verify explicitly the hypotheses of Brock's geometric theorem.
Definition~6.1 of~\cite[p.~193]{Brock2000} requires membership in
$S_+(\mathbb R^n)$ and continuous differentiability only on the open
set $\{0<v<\sup v\}$; Definition~6.2 on p.~194 requires local
symmetry in every direction. Theorem~6.1 on pp.~194--195 expressly
allows an unbounded positivity set and an infinite outer radius.
In particular, regularity at maximum points is not a hypothesis.
The global Sobolev assumption in the separate analytic
Theorem~6.2 is not needed for this geometric implication.

Set
\[
 v(x)=\frac{U(x)}{1+U(x)}\quad(x\ne0),\qquad v(0)=1.
\]
The two-sided bounds imply that $v$ is continuous on $\mathbb R^n$,
$\inf v=0$, $\sup v=1$, and $v(x)\to0$ as $|x|\to\infty$.
For a nonnegative function with infimum zero, membership in Brock's
$S_+(\mathbb R^n)$ means that every positive superlevel set has finite
Lebesgue measure; see~\cite[p.~163]{Brock2000}. Here, using $U(0)=+\infty$
only to describe superlevel sets,
\[
 \{v>t\}=\left\{U>\frac{t}{1-t}\right\}
 \subset B_{\left(C_0(1-t)/t\right)^{1/a}}
 \quad(0<t<1),
\]
and $\{v>t\}$ is empty for $t\ge1$. Thus $v\in S_+(\mathbb R^n)$.
Moreover,
\[
 D_v:=\{0<v<\sup v\}=\mathbb R^n\setminus\{0\}
 \quad\hbox{is open},\qquad v\in C^1(D_v).
\]
The origin belongs to $\{v=\sup v\}$, so neither differentiability
nor Lipschitz regularity there is required, including when $a<1$.

For $0<t<1$, the fiber components of $\{v>t\}$ are exactly those of
$\{U>s\}$ with $s=t/(1-t)$. At any paired endpoints $l,r$ with
$U_y(z,l)>0$,
\[
 \nabla v(z,l)=(1+s)^{-2}\nabla U(z,l),\qquad
 \nabla v(z,r)=(1+s)^{-2}\nabla U(z,r).
\]
Consequently~\eqref{eq:global-paired-derivatives}, including the
preceding exclusion of a critical opposite endpoint, verifies \cite[Definition~6.1]{Brock2000} for $v$ in every rotated coordinate system.
All hypotheses of~\cite[Theorem~6.1]{Brock2000} are therefore met.
Its conclusion transfers from $v$ to $U=v/(1-v)$ on $D_v$,
since the inverse transformation has strictly positive derivative.
We obtain the disjoint decomposition
\begin{equation}\label{eq:global-radial-pieces}
 \mathbb R^n\setminus\{0\}
 =\left(\bigcup_k\mathcal A_k\right)\cup S,\qquad
 \mathcal A_k=\{x:\rho_k<|x-z_k|<R_k\},
\end{equation}
where the open annuli are pairwise disjoint,
$0\leq\rho_k<R_k\leq\infty$, $U$ is radial and strictly decreasing
in radius on each annulus, and $\nabla U=0$ on $S$.

On one such annulus write $U(x)=\phi(|x-z|)$.
The radial weak equation makes the continuous flux
$P(r)=r^{n-1}|\phi'(r)|^{p-2}\phi'(r)$ locally absolutely continuous,
with
\begin{equation}\label{eq:global-radial-flux}
 P'(r)=-r^{n-1}\phi(r)^q<0,\qquad\phi'(r)<0.
\end{equation}
Its outer radius is infinite.  Otherwise choose $r_1\in(\rho,R)$.
For $r_1<r<R$ one has $0<\phi(r)\leq\phi(r_1)$ and
$P(r)\leq P(r_1)<0$.  The flux equation then gives a finite strictly
negative limit at $R$, so $\phi'(R^-)<0$.
If $0\in\partial B_R(z)$, approaching zero inside the annulus
would keep $U\leq\phi(r_1)$, contrary to the lower bound.
Otherwise $U$ is $C^1$ near every point of that sphere, with
nonzero gradient there.  By \eqref{eq:global-radial-pieces}, each
such point belongs to another open annulus, which must intersect
the original one.  This contradicts disjointness, proving $R=\infty$.

All annuli in \eqref{eq:global-radial-pieces} are therefore exteriors
of balls, and any two intersect.  There is at most one.
There is at least one because $\nabla U$ cannot vanish everywhere
in view of the strictly positive right side of the equation.
Write it as $\mathcal A=\{x:|x-z|>\rho\}$.
If $\rho>0$, the nonempty open set $B_\rho(z)\setminus\{0\}$ lies
in $S$, so its gradient vanishes.  Testing the weak equation with
a nonzero nonnegative smooth function supported there gives a
contradiction.  Thus $\rho=0$.
Finally $\mathbb R^n\setminus\{z\}=\mathcal A\subset
\mathbb R^n\setminus\{0\}$ forces $z=0$, proving the proposition.
\end{proof}

\section{Global radial symmetry without growth assumptions}
\label{sec:global-without-bounds}

We now remove the two-sided hypothesis from
Proposition~\ref{prop:two-sided-rigidity}. The argument uses two
scale-independent consequences of the upper-bound analysis.
We state and justify their arbitrary-scale versions explicitly.
Throughout this section, $U>0$ is a continuous weak solution of
$-\Delta_pU=U^{p^*-1}$ on $\mathbb R^n\setminus\{0\}$, and
\[
 m_U(r)=\min_{|x|=r}U(x),\qquad h_U(r)=r^a m_U(r).
\]
The Pohozaev constant of $U$, defined in
Lemma~\ref{lem:pohozaev}, is denoted by $\mathcal P(U)$.

\subsection{Two consequences of the blow-up analysis at arbitrary scales}

\begin{lemma}[Small critical minima force zero Pohozaev constant]
\label{lem:global-small-minimum}
If there is a sequence $t_j>0$ such that $h_U(t_j)\to0$, then
$\mathcal P(U)=0$. 
\end{lemma}

\begin{proof}
Set
\[
 w_j(y)=\frac{U(t_jy)}{m_U(t_j)},\qquad
 \varepsilon_j=h_U(t_j)^\sigma\longrightarrow0.
\]
Then
\begin{equation}\label{eq:global-scale-normalization}
 -\Delta_pw_j=\varepsilon_jw_j^q
 \quad\text{in }\mathbb R^n\setminus\{0\},\qquad
 \min_{\partial B_1}w_j=1.
\end{equation}
Extend $U$, and hence each $w_j$, $p$-superharmonically across zero.
For every fixed $R\ge2$, the minimum principle gives $w_j\ge1$
on $B_1$. On $B_{2R}\setminus\overline B_1$, compare with
\[
 \frac{|y|^{-b}-(2R)^{-b}}{1-(2R)^{-b}}.
\]
Nonnegativity on the outer sphere and the normalization on the
inner sphere give
\[
 w_j\ge c_R>0\quad\text{on }B_R.
\]
If $\mu_j$ is the Riesz measure of the extension, the weak Harnack
inequality and Lemma~\ref{lem:ub-mass}, using
$\inf_{B_{2R}}w_j\le1$, also give
\[
 \int_{B_R}w_j^\gamma\,dy\le C_R,\qquad
 \mu_j(B_R)\le C_R
 \quad\left(0<\gamma<\frac{n(p-1)}{n-p}\right).
\]
The measure estimate includes any atom at the origin.

These are exactly the inputs of the compactness argument
\eqref{eq:ub-basic-compactness}--\eqref{eq:ub-harmonic-limit}.
That argument applies here without a restriction on the sequence
$t_j$. Its only use of $t_j\downarrow0$ in the local problem was
to make the rescaled domains contain every fixed punctured ball;
the domains in \eqref{eq:global-scale-normalization} already have
this property. More explicitly, the proof of
\eqref{eq:ub-amplitude-cap}, using Lemma~\ref{lem:ub-tail} and the
positive lower comparison above, gives 
\[
 \sup_j\sup_K\varepsilon_j^a w_j<\infty, 
 \qquad K\Subset\mathbb R^n\setminus\{0\}. 
\]
The intrinsic core estimate \eqref{eq:ub-core-mass} then assigns
a locally uniform positive mass to each off-origin concentration
point. The Riesz measure bound makes their set $\Sigma$ locally
finite away from zero. Outside $\Sigma\cup\{0\}$, the intrinsic
potential is locally bounded; the bounded-potential Harnack
inequality and the preceding integral bound give uniform local
bounds for $w_j$. Interior regularity therefore gives, after
extraction,
\[
 w_j\longrightarrow W
 \quad\text{in }C^1_{\mathrm{loc}}
       (\mathbb R^n\setminus(\Sigma\cup\{0\})).
\]
None of these steps uses zero Pohozaev flux or the subsequent
classification of normalized harmonic limits.

Choose $R\in(1,2)$ so that $\partial B_R$ avoids $\Sigma$.
Local finiteness gives an open neighborhood of this sphere free
of concentration points. Consequently
\[
 \sup_j\left|\mathcal P_{\varepsilon_j}(w_j,R;0)\right|<\infty.
\]
The exact scaling identity is
\begin{equation}\label{eq:global-scale-pohozaev}
 \mathcal P(U)
 =m_U(t_j)^p t_j^{n-p}
       \mathcal P_{\varepsilon_j}(w_j,R;0)
 =h_U(t_j)^p\mathcal P_{\varepsilon_j}(w_j,R;0).
\end{equation}
Since $h_U(t_j)\to0$, this proves $\mathcal P(U)=0$.
\end{proof}

\begin{lemma}[Large critical values force small spherical minima]
\label{lem:global-large-value}
Suppose $x_j\ne0$, $r_j=|x_j|$, and
$K_j=r_j^aU(x_j)\to\infty$. Then, after extraction,
$h_U(r_j)\to0$. 
\end{lemma}

\begin{proof}
We verify the scale-independent application of the selection and
bubble-tail argument leading to \eqref{eq:ub-small-minimum}.
In $D_j=B_{r_j/4}(x_j)$ maximize $d^{n-p}U$, where
$d(x)=\operatorname{dist}(x,\partial D_j)$. Denote the selected
point, height, and boundary distance by $\xi_j,H_j,d_j$.
As in \eqref{eq:ub-selected-product},
\[
 H_j\ge U(x_j),\qquad
 H_jd_j^{n-p}\ge c r_j^{a(p-1)}K_j.
\]
Set
\[
 \delta_j=H_j^{-1/a},\qquad
 \Gamma_j=\frac{d_j}{8\delta_j},\qquad
 v_j(y)=\frac{U(\xi_j+\delta_jy)}{H_j}.
\]
The selected functions solve the critical equation, satisfy
$v_j(0)=1$ and $v_j\le2^{n-p}$ on $B_{4\Gamma_j}$, and
\[
 \Gamma_j^{n-p}=c d_j^{n-p}H_j^p\ge cK_j^p\longrightarrow\infty.
\]
Moreover, on each fixed compact set,
\[
 v_j(y)\le
 \left(\frac{d_j}{d_j-\delta_j|y|}\right)^{n-p}
 \longrightarrow1.
\]
Thus Zhang \cite[Theorem~1.1]{Zhang2026} and local compactness give
$v_j\to\mathcal U$ in $C^1_{\mathrm{loc}}(\mathbb R^n)$
along a subsequence.

Let
$\lambda_j=K_j^{-1/[2(n-p)]}$ and $S_j=\lambda_j\Gamma_j$.
Then $S_j\to\infty$ and $S_j=o(\Gamma_j)$.
Lemma~\ref{lem:ub-tail} gives, on the sphere of radius
$\delta_jS_j=d_j\lambda_j/8$ about $\xi_j$,
\[
 U\le CH_jS_j^{-b}
   =C\lambda_j^{-b}(H_jd_j^{n-p})^{-1/(p-1)}.
\]
This sphere lies in $D_j$. The superharmonic lower comparison
used in the previous lemma, at scale $r_j$, gives
$U\ge c\,m_U(r_j)$ there. Hence
\[
 h_U(r_j)
 \le C\lambda_j^{-b}K_j^{-1/(p-1)}
 =C K_j^{-1/[2(p-1)]}\longrightarrow0.
\]
This proves Lemma \ref{lem:global-large-value}.  
\end{proof}

\subsection{Strict negativity of the singular Pohozaev constant}

\begin{lemma}\label{lem:global-negative-pohozaev}
If the origin is a nonremovable singularity of $U$, then
$\mathcal P(U)<0$.
\end{lemma}

\begin{proof}
By Theorem~\ref{thm:main}, there exist constants $c_0,C_0,r_0>0$ such that
\[
 c_0|x|^{-a}\le U(x)\le C_0|x|^{-a}, 
 \qquad 0<|x|<r_0.
\]
Take $\rho_j\downarrow0$. Interior estimates and diagonal extraction
give
\[
 \rho_j^aU(\rho_j\,\cdot)\longrightarrow V
 \quad\text{in }C^1_{\mathrm{loc}}(\mathbb R^n\setminus\{0\}),
\]
where $V$ solves the same equation and satisfies the global bounds
\[
 c_0|x|^{-a}\le V(x)\le C_0|x|^{-a}, \qquad x\ne0.
\]
Proposition~\ref{prop:two-sided-rigidity} makes $V$ radial and
strictly decreasing. Critical scaling preserves the Pohozaev
flux, and $C^1$ convergence on a fixed sphere gives
\begin{equation}\label{eq:global-tangent-pohozaev}
 \mathcal P(U)=\mathcal P(V).
\end{equation}

Write $V(x)=v(|x|)$ and put
$\mu=\inf_{r>0}r^av(r)>0$.
Choose $s_k>0$ with $s_k^av(s_k)\to\mu$.
The radial solutions $s_k^aV(s_k\,\cdot)$ retain both global bounds.
A further $C^1_{\mathrm{loc}}$ limit is a radial solution
$Z(x)=z(|x|)$ such that
\[
 z(1)=\mu,\qquad r^az(r)\ge\mu\quad(r>0),\qquad
 \mathcal P(Z)=\mathcal P(V).
\]
For $f(r)=r^az(r)$, the minimum at $r=1$ gives
$f'(1)=0$ and $z'(1)=-a\mu<0$.
The radial equation is therefore nondegenerate near $1$, so that
$z$ is $C^2$ there. A direct calculation gives
\[
 0\le f''(1)
 =\frac{a^p\mu^{p-1}-\mu^q}
        {(p-1)(a\mu)^{p-2}}.
\]
It follows that $\mu^\sigma\le a^p$.
Evaluating the Pohozaev flux on $\partial B_1$ now yields
\begin{equation}\label{eq:global-negative-flux}
 \begin{aligned}
 \frac{\mathcal P(U)}{\omega_{n-1}}
 &=\frac{\mathcal P(Z)}{\omega_{n-1}}
   =-\frac{a^p\mu^p}{p}+\frac{\mu^{p^*}}{p^*}\\
 &\le -a^p\mu^p\left(\frac1p-\frac1{p^*}\right)
   =-\frac{a^p\mu^p}{n}<0.
 \end{aligned}
\end{equation}
This argument uses no classification of radial Fowler trajectories.
\end{proof}

\subsection{Proof of the global theorem}

\begin{proof}[Proof of Theorem~\ref{thm:global-radial}]
Nonremovability excludes $U\equiv0$, and the strong minimum
principle gives $U>0$. Lemma~\ref{lem:global-negative-pohozaev}
shows that $\mathcal P(U)<0$.

If $\inf_{r>0}h_U(r)=0$, choose $t_j>0$ with $h_U(t_j)\to0$.
Lemma~\ref{lem:global-small-minimum} would give $\mathcal P(U)=0$,
a contradiction. Therefore
\[
 c_*:=\inf_{r>0}h_U(r)>0,
 \qquad U(x)\ge c_*|x|^{-a}, \quad x\ne0.
\]
If $\sup_{x\ne0}|x|^aU(x)=\infty$, choose a sequence of such
unbounded critical values. Lemma~\ref{lem:global-large-value}
would produce radii with $h_U(r_j)\to0$, again a contradiction.
This proves the global upper bound and hence
\eqref{eq:global-automatic-bounds}.
Proposition~\ref{prop:two-sided-rigidity} now applies to $U$ itself
and gives the asserted radial symmetry and strict radial decrease.
\end{proof}

We now apply \cite[Proposition~2.11]{Franca2004Classification} to obtain the classification of radial singular solutions. To specify the normalization,
let $U(x)=\Phi(|x|)$ be the solution already obtained in
Section~\ref{sec:global-without-bounds}, and set
\begin{equation}\label{eq:fowler-transform}
 t=\log r,\qquad v(t)=r^a\Phi(r),\qquad
 z(t)=\bigl(-r^{a+1}\Phi'(r)\bigr)^{p-1}>0.
\end{equation}
The radial flux equation
\[
 \bigl(r^{n-1}|\Phi'|^{p-2}\Phi'\bigr)'
       =-r^{n-1}\Phi^q
\]
gives
\[
 \dot v=av-z^{1/(p-1)},\qquad \dot z=v^q-az.
\]
Its conserved quantity, in our normalization, is
\begin{equation}\label{eq:fowler-hamiltonian}
 H(v,z)=\frac{z^{p'}}{p'}-avz+\frac{v^{p^*}}{p^*},
 \qquad
 H(v(t),z(t))=\frac{\mathcal P(U)}{\omega_{n-1}}.
\end{equation}
Franca's variables are $(x,y)=(v,-z)$, the coefficient is $K\equiv1$, and the critical exponent denoted by $p^*$ is our $p^*-1$.

Proposition~2.11(C)--(D) of \cite{Franca2004Classification} gives
the negative-energy classification: the equilibrium at energy
$E_*=-a^n/n$, or the positive periodic orbit at energy $E\in(E_*,0)$,
up to time translation. Part~(B) identifies the zero-energy
homoclinic orbit with regular ground states. We need only the
negative-energy assertion. These interior orbits lie in
$v>0$, $z>0$, where the vector field is smooth for every $1<p<n$;
no differentiability on the phase-plane axes is required.

All hypotheses for the singular ground states in that reference
are satisfied here: radial symmetry and strict decrease have already
been proved, and \eqref{eq:global-automatic-bounds} implies
$\Phi(r)\to\infty$ at zero and $\Phi(r)\to0$ at infinity.
Moreover, Lemma~\ref{lem:global-negative-pohozaev} gives
$\mathcal P(U)<0$. Thus the cited result yields exactly the two
representations in Theorem~\ref{thm:global-radial}. Taking the minimum
of $v$ as time zero gives its stated normalization
$G(\alpha_E)=E$, $z=(a\alpha_E)^{p-1}$.
The energy and the phase modulo the minimal period are unique.
Conversely, the equilibrium and these periodic trajectories give
positive global singular solutions by the inverse of
\eqref{eq:fowler-transform}. The identity
\eqref{eq:fowler-hamiltonian} identifies their Pohozaev constants.
This completes the application of the classical radial classification
to Theorem~\ref{thm:global-radial}.

\section{Asymptotics to a Fowler solution}
\label{sec:asymptotic-symmetry} 
In this section, we prove Theorem~\ref{thm:asymptotic-radial-symmetry}. All constants
in this section may depend on the fixed solution and its Pohozaev
constant. For a function on the sphere, $\langle f\rangle$ denotes its
normalized spherical average. Recall that 
\[
 a=\frac{n-p}{p},\qquad p^*=\frac{np}{n-p},\qquad q=p^*-1,\qquad1<p<n.
\]

\subsection{Compactness and exclusion of the degenerate gradient set} 

Set
\begin{equation}\label{eq:asym-cylinder-definitions}
 t=\log r,\qquad v(t,\theta)=e^{at}u(e^t\theta),\qquad
 Dv=(v_t-av,\nabla_Sv),\qquad w=av-v_t.
\end{equation}
Here $\nabla_S$ is the spherical gradient, and $D$ also denotes
the linear operator $(\partial_t-a,\nabla_S)$ on test functions.
The change of variables, with the Euclidean test function
$r^{-a}\phi(\log r,\theta)$, gives the exact weak equation
\begin{equation}\label{eq:asym-cylinder-weak}
 \int\bigl(\A(Dv)\cdot D\phi-v^q\phi\bigr)\,dt\,d\theta=0,
 \qquad \A(\xi)=|\xi|^{p-2}\xi.
\end{equation}
Indeed, $p(a+1)=n$ and $ap^*=n$ cancel both radial weights.

Theorem~\ref{thm:main} and interior estimates on rescaled annuli give
\begin{equation}\label{eq:asym-cylinder-bounds}
 0<c\le v\le C,\qquad
 \sup_{s<T}\|v(s+\cdot,\cdot)\|_
 {C^{1,\alpha}([-2,2]\times\Sn)}<\infty
\end{equation}
for some $T<0$. Given any sequence $s_j\to-\infty$, diagonal
extraction gives a $C^1_{\mathrm{loc}}$ limit of
$v(s_j+\cdot,\cdot)$. In Euclidean variables it corresponds to
a positive global solution with the inherited two-sided critical
bounds. Its origin is nonremovable. By
Theorem~\ref{thm:global-radial}, this limit is radial and is either
the constant cylinder profile $a^a$ or a positive periodic Fowler
profile.

There are $T_1<0$ and $\kappa>0$ such that
\begin{equation}\label{eq:asym-nondegenerate}
 w(t,\theta)\ge\kappa, \qquad t<T_1,\ \theta\in\Sn.
\end{equation}
Otherwise, choose $s_j\to-\infty$ and $\theta_j\in\Sn$ with
$w(s_j,\theta_j)\le1/j$. Pass to a translated $C^1$ limit and
a convergent subsequence of $\theta_j$. Its radial Euclidean
representative $U$ is strictly decreasing, so the limiting value
of $w$ at time zero is $-U'(1)>0$, a contradiction.
This argument precedes every use of uniform ellipticity below.
Together with \eqref{eq:asym-cylinder-bounds}, it puts the scaled
gradient $Dv$ in a fixed compact subset of $\R^n\setminus\{0\}$. 
Let
\[
 E=\frac{\mathcal P(u)}{\omega_{n-1}}.
\]
The flux is conserved by Lemma~\ref{lem:pohozaev}, is invariant under
critical scaling, and passes to every translated limit by $C^1$
convergence on a sphere. The radial classification therefore gives
$E\in[E_*,0)$. If $E=E_*$, every translated limit is $a^a$.
A contradiction sequence and compactness give
\[
 v(t,\cdot)\longrightarrow a^a,\qquad
 v_t(t,\cdot)\longrightarrow0,\qquad
 \nabla_Sv(t,\cdot)\longrightarrow0
\]
uniformly as $t\to-\infty$. Both conclusions of the theorem follow
with $U_0(r)=a^ar^{-a}$.

Henceforth suppose $E_*<E<0$. Fix the normalized periodic profile
$V=V_E$, and write $T_0=T(E)$ and $g=aV-V'>0$. Every translated
limit is a phase translate of this same profile. Consequently,
for every fixed $L>0$,
\begin{equation}\label{eq:asym-local-phase-approximation}
 \lim_{s\to-\infty}\inf_{\tau\in\R/T_0\mathbb Z}
 \|v(s+\cdot,\cdot)-V(\cdot+\tau)\|_
 {C^1([-L,L]\times\Sn)}=0.
\end{equation}
If this failed, a sequence violating it would have a translated
limit of precisely the excluded form. No single phase for different
intervals is asserted in this step.

Put
\begin{equation}\label{eq:asym-mean-decomposition}
 m(t)=\langle v(t,\cdot)\rangle,\qquad
 \eta=v-m,\qquad s(t)=am-m'=\langle w\rangle.
\end{equation}
Then $s\ge\kappa$, and
\begin{equation}\label{eq:asym-angular-small}
 \|\eta(t,\cdot)\|_{L^\infty}
 +\|\eta_t(t,\cdot)\|_{L^\infty}
 +\|\nabla_S\eta(t,\cdot)\|_{L^\infty}\longrightarrow0
 \quad(t\to-\infty).
\end{equation}

\subsection{A strict angular coercivity estimate}

For compactly supported functions on $\mathcal C=\R\times\Sn$,
define the quadratic form at the radial model by
\begin{equation}\label{eq:asym-model-form}
 Q_V[\phi]=\int_{\mathcal C}
 \left\{(p-1)g^{p-2}(\phi_t-a\phi)^2
       +g^{p-2}|\nabla_S\phi|^2-qV^{q-1}\phi^2\right\}.
\end{equation}

\begin{lemma}\label{lem:asym-model-gap}
There is $c_0>0$ such that
\[
 Q_V[\phi]\ge c_0\|\phi\|_{H^1(\mathcal C)}^2
 \quad\text{whenever }\langle\phi(t,\cdot)\rangle=0
 \text{ for every }t.
\]
The same constant works for every phase translate of $V$.
\end{lemma}

\begin{proof}
Write
\[
 A=(p-1)g^{p-2},\qquad B=aA'+a^2A-qV^{q-1}.
\]
The positive function $Y=e^{-t}g$ satisfies
\begin{equation}\label{eq:asym-translation-jacobi}
 -(AY')'+\bigl(B+(n-1)g^{p-2}\bigr)Y=0.
\end{equation}
One way to verify this identity is to differentiate the equation
for $U(r)=r^{-a}V(\log r)$ under Euclidean translations:
$r^a\partial_{x_i}U=-e^{-t}g\,\theta_i$, and
$-\Delta_S\theta_i=(n-1)\theta_i$. The differentiation is legitimate
on every compact annulus because $U'(r)<0$. Alternatively,
substitute $z=g^{p-1}$, $z'=V^q-az$ and $V'=aV-g$
directly into \eqref{eq:asym-translation-jacobi}.

For a scalar compactly supported $\psi=gf$, integration by parts
in \eqref{eq:asym-translation-jacobi} gives
\begin{align}
 &\int_{\R}\left\{A\psi'^2+
           \bigl(B+(n-1)g^{p-2}\bigr)\psi^2\right\}\,dt
 \notag\\
 &\qquad=\int_{\R}AY^2
          \left|\left(\frac{\psi}{Y}\right)'\right|^2\,dt
       =(p-1)\int_{\R}g^p(f'+f)^2\,dt
 \notag\\
 &\qquad\ge(p-1)(\min g)^p
            \int_{\R}(f'^2+f^2)\,dt.                 \label{eq:asym-ground-state-gap}
\end{align}
The last inequality first bounds the weight from below and then uses
$\int 2ff'=0$; no weighted cross term is discarded.
Periodic positivity and smoothness of $g$ imply that the last
expression controls $\|\psi\|_{H^1(\R)}^2$.

For a slice-mean-zero $\phi$, the spherical Poincar\'e inequality
with first eigenvalue $n-1$ shows that
\[
 Q_V[\phi]\ge c\int_{\mathcal C}(\phi_t^2+\phi^2).
\]
On the other hand, bounded coefficients and the positive principal
coefficients in \eqref{eq:asym-model-form} give
\[
 Q_V[\phi]\ge c'\int_{\mathcal C}
       (|\phi_t|^2+|\nabla_S\phi|^2)-C'\int_{\mathcal C}\phi^2.
\]
A suitable convex combination proves the full $H^1$ bound.
All bounds on the periodic coefficients are unchanged by a phase
translation.
\end{proof}

We now transfer this estimate to the actual solution. Define
\begin{equation}\label{eq:asym-secant-coefficients}
 \mathbb B(t,\theta)
   =\int_0^1 D\A(Dm+\rho D\eta)\,d\rho,\qquad
 c_*(t,\theta)=q\int_0^1(m+\rho\eta)^{q-1}\,d\rho.
\end{equation}
Every point on each gradient segment has time component
\[
 (Dm+\rho D\eta)_t=-(1-\rho)s-\rho w\le-\kappa.
\]
All these segments lie in a fixed bounded set away from zero.
Thus
\begin{equation}\label{eq:asym-uniform-ellipticity}
 \lambda|\xi|^2\le\mathbb B\xi\cdot\xi\le\Lambda|\xi|^2,
 \qquad |c_*|\le C,
\end{equation}
for every $1<p<n$. Indeed,
\[
 D\A(\xi)=|\xi|^{p-2}
 \left(I+(p-2)\frac{\xi\otimes\xi}{|\xi|^2}\right)
 \quad(\xi\ne0).
\]
Also, $m+\rho\eta\ge c>0$, so $q<1$ causes no difficulty.

For every compactly supported slice-mean-zero test function $\phi$,
\[
 \int\bigl(\A(Dm)\cdot D\phi-m^q\phi\bigr)=0,
\]
because all coefficients depend only on $t$. This identity does
not assert that $m$ is a radial solution. Subtracting it from
\eqref{eq:asym-cylinder-weak} gives the exact projected equation
\begin{equation}\label{eq:asym-projected-equation}
 \int\left(\mathbb B D\eta\cdot D\phi-c_*\eta\phi\right)=0.
\end{equation}
Let
\[
 Q[\phi]=\int\left(\mathbb B D\phi\cdot D\phi-c_*\phi^2\right).
\]

\begin{lemma}\label{lem:asym-tail-gap}
After decreasing $T_1$, there is $c_1>0$ such that
\[
 Q[\phi]\ge c_1\|\phi\|_{H^1}^2
\]
for every compactly supported slice-mean-zero $\phi$ in
$(-\infty,T_1)\times\Sn$.
\end{lemma}

\begin{proof}
Choose a smooth partition depending only on time with
\[
 \sum_j\chi_j^2=1,\qquad
 \operatorname{diam}(\supp\chi_j)\le C L,\qquad
 \sum_j|\chi_j'|^2\le C L^{-2}.
\]
Its overlap is bounded independently of $L$. Direct expansion,
using $D(\chi_j\phi)=\chi_jD\phi+\chi_j'\phi e_t$, yields
\begin{equation}\label{eq:asym-ims}
 \sum_jQ[\chi_j\phi]
 =Q[\phi]+\int\mathbb B_{tt}
                 \sum_j|\chi_j'|^2\phi^2.
\end{equation}
Fix $L$ sufficiently large that this last error is smaller than
one quarter of the model coercivity constant. Then use
\eqref{eq:asym-local-phase-approximation} on intervals of length
$O(L)$. On each sufficiently negative interval, the coefficients
in \eqref{eq:asym-secant-coefficients} are uniformly as close as
desired to those of one phase of $Q_V$. Smoothness of $\A$ on the
nondegenerate segments and positivity of $m,v$ justify this coefficient
comparison. It follows that
\[
 Q[\chi_j\phi]\ge\frac{c_0}{2}\|\chi_j\phi\|_{H^1}^2.
\]
Different intervals may use different phases. Only the
phase-independent constant in Lemma~\ref{lem:asym-model-gap} is needed.
Move the finite endpoint of the tail left by a fixed multiple of
$L$ so that every relevant partition interval lies in the region of
comparison. Summing the preceding inequality and using
\eqref{eq:asym-ims} proves the assertion, since
$\sum_j\|\chi_j\phi\|_{H^1}^2\ge\|\phi\|_{H^1}^2$.
\end{proof}

\subsection{Integrability of the angular error}

A time cutoff preserves the slice-mean-zero condition.
Testing \eqref{eq:asym-projected-equation} with $\chi^2\eta$
and expanding gives the exact identity
\begin{equation}\label{eq:asym-cutoff-identity}
 Q[\chi\eta]=\int\mathbb B_{tt}|\chi'|^2\eta^2.
\end{equation}
These tests are admissible by approximation on compact strips;
$\eta$ is locally $C^1$, and all coefficients in this identity are
bounded.

Fix $T_2<T_1$ and a right cutoff which is one on $(-\infty,T_2]$
and vanishes before $T_1$. Multiply it by a left cutoff which
vanishes on $(-\infty,-2R]$, is one on $[-R,\infty)$, and has
derivative bounded by $C/R$. For large $R$ the two transition
regions are disjoint. The left contribution to
\eqref{eq:asym-cutoff-identity} is $O(R^{-1})$ by boundedness of
$\eta$; the right contribution is a fixed finite constant.
Lemma~\ref{lem:asym-tail-gap} and Fatou's lemma therefore give
\begin{equation}\label{eq:asym-angular-integrable}
 \int_{-\infty}^{T_2}e(t)\,dt<\infty,\qquad
 e(t)=\left\langle\eta^2+\eta_t^2+|\nabla_S\eta|^2\right\rangle.
\end{equation}
Moreover, $e(t)\to0$ by \eqref{eq:asym-angular-small}, independently
of the integrability argument.  

\subsection{The averaged flux and the direct energy error}

Define
\begin{equation}\label{eq:asym-averaged-flux}
 S=(w^2+|\nabla_Sv|^2)^{1/2},\qquad
 J=S^{p-2}w,\qquad Z(t)=\langle J(t,\cdot)\rangle.
\end{equation}
Equation \eqref{eq:asym-cylinder-weak} says in distributions that
\[
 J_t+aJ-\diver_S(S^{p-2}\nabla_Sv)=v^q.
\]
Averaging gives
\begin{equation}\label{eq:asym-exact-average}
 Z'=\langle v^q\rangle-aZ.
\end{equation}
In particular, $Z$ is $C^1$ since the right-hand side is continuous.
The variables $m,s,Z$ all stay in fixed compact subintervals of
$(0,\infty)$ on the tail.

Write $d=w-s=a\eta-\eta_t$, so $\langle d\rangle=0$.
Taylor expansion at $(s,0)$, on the nondegenerate segments already
specified, gives
\[
 ((s+d)^2+|\nabla_S\eta|^2)^{(p-2)/2}(s+d)
 =s^{p-1}+(p-1)s^{p-2}d
       +O(d^2+|\nabla_S\eta|^2).
\]
The Hessians here are uniformly bounded even when $1<p<2$,
because the radial component stays at least $\kappa$.
Likewise, expansion of $(m+\eta)^q$ is uniform even when $q<1$,
because $m+\rho\eta\ge c$. Averaging cancels the linear terms:
\begin{equation}\label{eq:asym-quadratic-average}
 Z=s^{p-1}+O(e(t)),\qquad
 \langle v^q\rangle=m^q+O(e(t)).
\end{equation}
Since the inverse power is smooth on the attained positive interval,
$X=(m,Z)$ satisfies
\begin{equation}\label{eq:asym-perturbed-system}
 X'=F(X)+R(t),\qquad
 F(m,Z)=(am-Z^{1/(p-1)},\,m^q-aZ),\qquad
 |R(t)|\le C e(t).
\end{equation}

An additional estimate is essential: integrability of $R$ alone
does not settle phase convergence. Substitution of
\eqref{eq:asym-cylinder-definitions} into the conserved Pohozaev
flux gives
\begin{equation}\label{eq:asym-cylinder-pohozaev}
 E=\left\langle
       -avJ+wJ-\frac{S^p}{p}+\frac{v^{p^*}}{p^*}
    \right\rangle.
\end{equation}
For the Hamiltonian
\[
 H(m,Z)=\frac{Z^{p'}}{p'}-amZ+\frac{m^{p^*}}{p^*},
\]
the value of the integrand at $(v,w,\nabla_Sv)=(m,s,0)$
is $H(m,s^{p-1})$. Its first variations in $v$ and $w$ have
zero average, and its first derivative in the angular-gradient
variable is zero at this radial point. All second derivatives
are uniformly bounded on the Taylor segments. Consequently,
\[
 E=H(m,s^{p-1})+O(e(t)).
\]
Together with \eqref{eq:asym-quadratic-average}, this proves
\begin{equation}\label{eq:asym-direct-energy-error}
 |H(X(t))-E|\le C e(t).
\end{equation}
Thus the energy deviation itself, not merely its derivative, is
integrable on the negative half-line.

\subsection{Convergence of the phase and completion of the proof}

Choose a closed interval $I\Subset(E_*,0)$ whose interior contains
$E$. The energy ovals with energy in $I$ form a compact annulus
strictly inside the positive quadrant. The vector field $F$ is smooth
there for every $1<p<n$; no regularity at either coordinate axis
is required.

For $h\in I$, let $\alpha(h)\in(0,a^a)$ be the smaller root of
$G(\alpha)=h$ and put
\[
 P(h)=(\alpha(h),(a\alpha(h))^{p-1}).
\]
Since $G'(\alpha(h))\ne0$, this is a smooth transverse section:
$H(P(h))=h$, whereas $dH(F)=0$ and $F$ does not vanish on these
ovals. Let $\Phi_\ell$ be the flow of $F$. The minimal return
period $T(h)$ is smooth, after shrinking $I$ if necessary.
Indeed, the return to a short transverse section near $P(E)$ is
transverse, so the implicit function theorem gives a smooth return
time; the simple closed ovals identify it with their minimal period.

The map
\begin{equation}\label{eq:asym-phase-coordinate}
 \Psi(h,\vartheta)=
 \Phi_{\vartheta T(h)/T_0}(P(h)),\qquad
 \vartheta\in\R/T_0\mathbb Z,
\end{equation}
is a smooth coordinate on this annulus. It is one-to-one because
$h$ is the energy and each oval is traversed exactly once as
$\vartheta$ runs around the circle; its derivative is nonsingular
because the energy and flow directions are transverse.
Write its inverse as $(H,\Theta)$. The differential of its
circle-valued angle satisfies
\[
 d\Theta(F)=\frac{T_0}{T(H)}.
\]
By \eqref{eq:asym-direct-energy-error} and $e(t)\to0$, $X(t)$ lies
in this annulus for all sufficiently negative $t$. Lift
$\Theta(X(t))$ continuously to a real-valued $C^1$ function
$\vartheta(t)$ on that half-line. Equations
\eqref{eq:asym-perturbed-system} and
\eqref{eq:asym-direct-energy-error} give
\[
 \left|\vartheta'(t)-1\right|
 \le C\bigl(|H(X(t))-E|+|R(t)|\bigr)
 \le C e(t).
\]
By \eqref{eq:asym-angular-integrable}, there is $\tau\in\R$ with
\begin{equation}\label{eq:asym-phase-limit}
 \vartheta(t)-t\longrightarrow\tau
 \qquad(t\to-\infty).
\end{equation}
Use \eqref{eq:asym-phase-coordinate} on its compact circle variable
to reconstruct the state. Since $H(X(t))\to E$, it follows that
\[
 X(t)-(V_E(t+\tau),Z_E(t+\tau))\longrightarrow0.
\]
This reconstruction uses a circle coordinate, so it introduces
no factor of the unbounded lifted phase multiplying the energy
error.

By \eqref{eq:asym-perturbed-system} and $e(t)\to0$,
\[
 m'(t)-V_E'(t+\tau)\longrightarrow0.
\]
Combining these limits with \eqref{eq:asym-angular-small} yields
\begin{equation}\label{eq:asym-cylinder-final}
 \sup_{\theta\in\Sn}\left(
 |v(t,\theta)-V_E(t+\tau)|
 +|v_t(t,\theta)-V_E'(t+\tau)|
 +|\nabla_Sv(t,\theta)|
 \right)\longrightarrow0.
\end{equation}
Define
\[
 U_0(r)=r^{-a}V_E(\log r+\tau).
\]
Then $U_0$ is a global singular solution.
Its profile has a positive minimum; thus
\eqref{eq:asym-cylinder-final} gives
\eqref{eq:asymptotic-fixed-profile}. 
This completes the proof of
Theorem~\ref{thm:asymptotic-radial-symmetry}.

\bigskip

\noindent\textbf{Acknowledgements} Part of this work was completed while the author was visiting The Chinese University of Hong Kong. The author would like to thank Professor Juncheng Wei for his support and encouragement, as well as the Department of Mathematics and the Institute of Mathematical Sciences for their hospitality. The research of the author is supported by NSFC 12671139 and 12301140. The author used AI-assisted tools, principally ChatGPT/Codex, for brainstorming and manuscript editing. AI output was not treated as a mathematical authority. The author independently worked and verified all theorem statements, proofs, calculations, constants, citations, and final prose, and take full responsibility for the contents of the paper.

\medskip 

\noindent Hui Yang

\noindent School of Mathematical Sciences \\
Shanghai Jiao Tong University, Shanghai 200240, China \\[1mm]
Email: \textsf{hui-yang@sjtu.edu.cn}

\end{document}